\documentclass[12pt,a4paper,reqno]{amsart}
\usepackage{ifpdf}
\ifpdf
\usepackage{cmap}
\fi

\usepackage[T1]{fontenc}
\usepackage[utf8x]{inputenc}
\usepackage{amsthm}
\usepackage{amssymb}
\usepackage{amsmath}
\usepackage{bm}
\usepackage{soul}
\usepackage{sidecap}

\usepackage{fourier}

\usepackage[usenames,svgnames]{xcolor}
\usepackage{float}
\usepackage{graphicx}
\ifpdf
\usepackage[colorlinks=true,citecolor=DarkGreen,linktoc=page,unicode=true,pdftex]{hyperref}
\DeclareGraphicsExtensions{.pdf}
\else
\usepackage[colorlinks=true,citecolor=DarkGreen,linktoc=page,unicode=true
]{hyperref}
\DeclareGraphicsExtensions{.eps}
\fi
\usepackage{extarrows}
\usepackage{colonequals}

\newcommand{\red}{\color[HTML]{E00500}}
\newcommand{\bla}{\color{Black}}
\newcommand{\navy}{\color{Navy}}

\newcommand{\gray}{\color[HTML]{506070}}
\newcommand{\green}{\color[HTML]{509000}}

\newcommand{\Arg}{\operatorname{Arg}}
\newcommand{\sign}{\operatorname{sign}}

\renewcommand{\le}{\leqslant}
\renewcommand{\ge}{\geqslant}

\renewcommand{\Re}{\operatorname{Re}}
\renewcommand{\Im}{\operatorname{Im}}

\newtheorem{conjecture}{Conjecture}

\newtheorem{theorem}{Theorem}

\newtheorem{corollary}[theorem]{Corollary}

\numberwithin{paragraph}{section}
\author{Alexander Dyachenko}
\thanks{This work was financially supported by the European Research Council under the European
    Union's Seventh Framework Programme (FP7/2007--2013)/ERC grant agreement no. 259173.}
\title[On certain class of entire functions]{On certain class of entire functions and a
    conjecture by Alan Sokal}
\email{\href{mailto:dyachenk@math.tu-berlin.de}{dyachenk@math.tu-berlin.de}\\\href{mailto:diachenko@sfedu.ru}{diachenko@sfedu.ru}}
\address{TU-Berlin\\MA 4-2\\Straße des 17. Juni 136\\10623 Berlin\\Germany}

\begin{document}

\begin{abstract}
    In this paper we determine a class of entire functions using conditions on their odd and
    even parts. Further it is shown that the zeros of members of this class are localized in a
    very special way. This result allows us to treat a particular case of the conjecture by
    A.~Sokal (to show that all zeros of the function
    $\displaystyle {\sum}_{k=0}^\infty q^{\frac{k(k-1)}2} z^{k}\big/k!$, where $q\in\mathbb C$,
    $0<|q|\le 1$, are simple and distinct in absolute value).
\end{abstract}

\maketitle
\tableofcontents
\section{Introduction}
For a meromorphic function $F$ of the complex variable $z$ consider its decomposition into odd
and even parts such that
\begin{equation}\label{eq:f_even_odd}
F(z)=f(z^2) + z g(z^2).
\end{equation}
It is a well-known fact that for
stability\footnote{The polynomial is called \emph{stable} if all of its roots have negative real
    part.}
of a real\footnote{The function is \emph{real} whenever it is real on the real line.}
polynomial~$F$ it is necessary and sufficient that the functions $f(z)$ and $g(z)$ only have
simple negative
interlacing\footnote{The zeros of two functions are called \emph{interlacing} if in between two
    consecutive zeros of the first function there is a zero of the second function and
    \emph{vice versa.}} zeros and the same sign at the origin. This condition can be restated as
follows: the polynomials are coprime and their ratio $\phantom{\bigg.}\frac{f(z)}{g(z)}$ has
only negative zeros and poles and is positive at the origin. This correspondence (expressed as
conditions on the Hurwitz matrix) is at the heart of the Routh-Hurwitz theory (see,
\emph{e.g.}~\cite[Ch.~XV]{Gantmakher},~\cite{Tyaglov,BarkovskyTyaglov,AkhKrein,ChebMei}). With a
proper extension of the notion of stability, this criterion can be generalized into entire
(see~\cite{AkhKrein,ChebMei}), rational (see~\cite{BarkovskyTyaglov}) and further toward
meromorphic functions. Furthermore, if $F(z)$ is a polynomial and we allow the ratio
$\phantom{\bigg.}\frac{f(z)}{g(z)}$ to have positive zeros and poles, then we will obtain the
``generalized Hurwitz'' polynomials (in terms of the work~\cite{Tyaglov}). In the same
publication (see \cite[Subsection~4.6]{Tyaglov}) the author describes ``strange polynomials''
(which are related to stable polynomials) with a very specific behaviour. In the present paper
the nature of this phenomenon is explained (Theorem~\ref{th:ca1}, the case~\ref{distr3}).

In fact, here we consider a more general question. How are the zeros of a (complex) entire
function $F(z)$ distributed if the ratio~$\phantom{\bigg.}\frac{f(z)}{g(z)}$ has only negative
zeros and positive poles? Or even more generally, how are the zeros of a meromorphic function
$F(z)$ distributed if both functions $f(z)$ and $g(-z)$ have only negative zeros and positive
poles? The obtained solution to the question allows us to treat a special case of a conjecture,
which was stated by Alan Sokal (see Section~\ref{sec:conj-alan-sokal}). This research originally
began with considering that conjecture.

The author is grateful to the people who made helpful comments concerning this publication,
especially to the members (former and current) of his working group at the TU-Berlin. I also
thank the colleagues from Potsdam and Ufa for useful discussions.

\section{Main result}
At first, we state our result for entire functions $F(z)$ of genus one at the most. This is
equivalent to the condition that the genus of the functions $f(z)$ and $g(z)$ defined
by~\eqref{eq:f_even_odd} is equal to zero. Then we discuss further extensions and alternative
conditions.

\subsection{Entire functions of genus one at the most}
Let entire functions $f(z)$ and $g(z)$ of genus 0 satisfy $f(z)=0\implies z<0$ and
$g(z)=0 \implies z\ge 0$, that is
\[
f(z) = f_0 \prod_{\mu=0}^{\omega_1}\left(1+\frac{z}{b_\mu}\right)
\quad\text{and}\quad
g(z) = g_0 z^j\prod_{\nu=0}^{\omega_2}\left(1-\frac{z}{a_\nu}\right),
\]
where $f_0\ne 0$, $g_0\ne 0$, $j$ is a nonnegative integer, $b_\mu,a_\nu>0$ for all
$\mu=0,\dots,\omega_1$, $\nu =0,\dots,\omega_2$, $0\le\omega_1\le \infty$,
$0 \le \omega_2 \le \infty$, and
$\sum_{\mu=0}^{\omega_1} b_\mu^{-1} + \sum_{\nu=0}^{\omega_2} a_\nu^{-1} < +\infty$.

We also allow the situation when one of the functions $f(z)$, $z^jg(z)$ is constant (not both of
them). However, we do not highlight this case specially, since the corresponding notation can be
abusing. If both of these functions are constants, then the zeros of~$F(z)=f_0+z^{2j+1}g_0$
satisfy $z^{2j+1} =-\frac{f_0}{g_0} = \frac 1c \frac{|f_0|}{|g_0|}$; so Theorem~\ref{th:ca1} in
such case holds only for~$j=0$.

According to~\eqref{eq:f_even_odd}, the function $F(z)$ is equal to $f(z^2)+z g(z^2)$, so that
$f(z^2)$ and $zg(z^2)$ are its even and odd parts, respectively.

Denote $c\colonequals-\frac{g_0|f_0|}{f_0|g_0|}$ and
\begin{equation}\label{eq:R_form}
    R(z)\colonequals-c\frac{f(z^2)}{z g(z^2)}
    =\frac{1}{z^{2j+1}}\cdot
    \frac{|f_0|\prod_{\mu=0}^{\omega_1}\left(1+\frac{z^2}{b_\mu}\right)}
    {|g_0|\prod_{\nu=0}^{\omega_2}\left(1-\frac{z^2}{a_\nu}\right)},
\end{equation}
here the products on $\mu$ and $\nu$ are convergent. Using the notation~``$\arg$'' for the
multivalued argument function and~``$\Arg$'' for the principal branch of argument,
$-\pi<\Arg z\le\pi$ for any~$z$, we have
\[|c|=1 \quad\text{and}\quad \arg c = \pi + \arg g_0 - \arg f_0.\]
Furthermore,
\begin{equation}\label{eq:arg_R_real}
    \Arg R(x) \in \{0,\pi\}
    \quad\text{and}\quad
    \Arg R(ix) = \left\{-\tfrac{\pi}{2},\tfrac{\pi}{2}\right\}
    \quad\text{whenever}\quad x\in\mathbb{R}
\end{equation}

For the (open) quadrants of the complex plane we use the notation
\begin{align*}
    &Q_1\colonequals\{z\in\mathbb{C}:\Re z>0,\Im z>0\},&
    &Q_2\colonequals\{z\in\mathbb{C}:\Re z<0,\Im z>0\},&\\
    &Q_3\colonequals\{z\in\mathbb{C}:\Re z<0,\Im z<0\},&
    &Q_4\colonequals\{z\in\mathbb{C}:\Re z>0,\Im z<0\},&
\end{align*}
so that they are numbered in a counterclockwise direction. Observe that $\Im z^2> 0$ is the
same as $z\in Q_1\cup Q_3$, and $\Im z^2< 0$ is the same as $z\in Q_2\cup Q_4$. We also denote
the set of all zeros of the function~$F(z)$ under consideration via~$\{z_k\}_{k=1}^\omega$,
$1\le\omega\le\infty$, where each multiple zero is taken as many times as its multiplicity.

\begin{theorem}\label{th:ca1}\navy
    The zeros $\{z_k\}_{k=1}^\omega$ are distributed as follows:
    \begin{enumerate}
    \item\label{distr1} If\/~$\Im c^2>0$, then all zeros are simple and satisfying
        $0<|z_1|<|z_2|<\dots$; $\Im z_k^2\ne 0$~for every natural~$k$, if~$z_k\in Q_l$,
        $l=1,2,3$, then $z_{k+1}\in Q_{l+1}$ and if $z_k\in Q_4$ then $z_{k+1}\in Q_1$.
        Moreover, 
        $(-1)^j\Re c\Re z_1>0$ and $\Re z_1\Im z_1<0$.
    \item\label{distr2} If\/~$\Im c^2<0$, then all zeros are simple and satisfying
        $0<|z_1|<|z_2|<\dots$; $\Im z_k^2\ne 0$~for every natural~$k$, if~$z_k\in Q_l$,
        $l=2,3,4$, then $z_{k+1}\in Q_{l-1}$ and if $z_k\in Q_1$ then $z_{k+1}\in Q_4$.
        Moreover, 
        $(-1)^j\Im c\Im z_1 < 0$ and $\Re z_1\Im z_1>0$.
    \item\label{distr3} If\/~$\Im c=0$, then $0<|z_1|\le|z_2|<|z_3|\le|z_4|<|z_5|\dots$, there are
        no purely imaginary zeros, real zeros can be simple or double, other zeros are simple.
        Moreover, for each natural $k$ we have
        \begin{gather*}
            |z_{2k-1}|=|z_{2k}| \implies z_{2k-1}=\overline{z}_{2k},\quad
            |z_{2k-1}|<|z_{2k}| \implies \Arg z_{2k-1}=\Arg{z}_{2k}\in\{0,\pi\}\\
            \sign\Re z_{2k}=-\sign\Re z_{2k+1},\quad
            (-1)^jc\Re z_1 > 0 \quad\text{and}\quad 
            \Re z_1\Im z_1 \le 0;
        \end{gather*}
    \item\label{distr4} If\/~$\Re c=0$, then $0<|z_1|<|z_2|\le|z_3|<|z_4|\le|z_5|\dots$, there are
        no real zeros, purely imaginary zeros can be simple or double, other zeros are simple.
        Moreover, for each natural $k$ we have
        \begin{gather*}
            |z_{2k}|=|z_{2k+1}| \implies z_{2k}=-\overline{z}_{2k+1},\quad
            |z_{2k}|<|z_{2k+1}| \implies
                    \Arg z_{2k}=\Arg{z}_{2k+1}\in\left\{-\frac\pi2,\frac\pi2\right\}\\
            \sign\Im z_{2k-1}=-\sign\Im z_{2k},\quad
            (-1)^j\Im c\Im z_1 < 0 \quad\text{and}\quad \Re z_1 = 0.
        \end{gather*}
    \end{enumerate}
\end{theorem}

The distribution~\ref{distr3} reflects that the function~$F(z)$ is real up to a constant for
real~$c$, and the distribution~\ref{distr4} reflects that $F(iz)$ is real up to a constant for
purely imaginary~$c$.

Note that, in particular, the distributions~\ref{distr1}--\ref{distr4} are ``almost uniform''.
That is, if zeros of a function are distributed as in~\ref{distr1} or~\ref{distr2}, then for any
$r>0$ the number of zeros in $\{z\in Q_k:|z|<r\}$ can differ from the number of zeros in
$\{z\in Q_j:|z|<r\}$ at most by $1$ (here $k,j=1,\dots,4$). The cases~\ref{distr3}
and~\ref{distr4} are the ``degenerated'' cases \ref{distr1} and~\ref{distr2}, with possible
ingress of some zeros onto the real or imaginary axes.

\subsection{Some straightforward corollaries}
The standard test on negativity of zeros of an entire function is Grommer's criterion (see,
\emph{e.g.}~\cite{AkhKrein,ChebMei}). Applying this criterion immediately gives the following
conditions on coefficients of the expansion into the power series.
\begin{corollary}\label{cr:ca3}\navy
    Let an entire function $F(z)$ be of genus no more than one.  Suppose that there
    are nonzero complex numbers $\alpha$, $\beta$ and a nonnegative integer $j$, such that the
    function $F(z)$ has the following expansion into power series
    \begin{equation}\label{eq:pow_series}
        F(z)=\alpha \sum_{k=0}^\infty f_{2k} z^{2k} + \beta \sum_{l=0}^\infty (-1)^l f_{2l+1} z^{2l+2j+1},
    \end{equation}
    and the matrices
    \begin{equation}
        \label{eq:grommer}
        \begin{pmatrix}
            f_0& f_2& f_4& f_6& f_8&\hdots\\
            0  & f_2&2f_4&3f_6&4f_8&\hdots\\
            0  & f_0& f_2& f_4& f_6&\hdots\\
            0  & 0  & f_2&2f_4&3f_6&\hdots\\
            0  & 0  & f_0& f_2& f_4&\hdots\\
            \vdots&\vdots&\vdots&\vdots&\vdots&\ddots
        \end{pmatrix}
        \quad\text{and}\quad
        \begin{pmatrix}
            f_1& f_3& f_5& f_7& f_9&\hdots\\
            0  & f_3&2f_5&3f_7&4f_9&\hdots\\
            0  & f_1& f_3& f_5& f_7&\hdots\\
            0  & 0  & f_3&2f_5&3f_7&\hdots\\
            0  & 0  & f_1& f_3& f_5&\hdots\\
            \vdots&\vdots&\vdots&\vdots&\vdots&\ddots
        \end{pmatrix}
    \end{equation}
    have its leading principle minors positive (or positive up to some order and then all equal
    to zero, $\exists n>1: f_n\ne 0$). Then the function $F(z)$ satisfies Theorem~\ref{th:ca1} with
    $c=-\frac{\beta|\alpha|}{\alpha|\beta|}$.

    Conversely, if the function $F(z)$ satisfies Theorem~\ref{th:ca1} then for $\alpha=f_0$ and
    $\beta=g_0$ it has the form~\eqref{eq:pow_series} and the matrices~\eqref{eq:grommer} have
    their leading principle minors positive (or positive up to some order and then all equal to
    zero, $\exists n>1: f_n\ne 0$).
\end{corollary}

The number of minors to be checked can be reduced if we take advantage of the results by
Li\'enard and Chipart. The conditions on the odd or on the even minors can be replaced with
conditions on the coefficients. Furthermore, in (and only in) the case of a polynomial $F(z)$,
the number of conditions in Corollary~\ref{cr:ca3} is finite (see,
\emph{e.g.}~\cite{Gantmakher,Tyaglov}).

Observe that if $\mu^4=-1$ (\emph{i.e.} $\mu$ is a primitive $8$th root of unity) and
$\beta=\alpha\overline\mu$, then we have the identity
\begin{multline*}
    F(z)=\alpha \sum_{k=0}^\infty f_{2k} z^{2k}
         +\alpha \overline\mu \sum_{l=0}^\infty (-1)^k f_{2k} z^{2k+1}
     =\alpha \left(\sum_{k=0}^\infty f_{2k} 1^{\frac{k(k-1)}{2}} z^{2k}
         +\sum_{l=0}^\infty f_{2l+1} (-1)^{l^2}1^l\mu^{-1} z^{2l+1}\right)\\
     =\alpha \left(\sum_{k=0}^\infty f_{2k} \mu^{4k(k-1)} z^{2k}
         +\sum_{l=0}^\infty f_{2l+1} \mu^{4l^2-1} z^{2l+1}\right)
     =\alpha \sum_{k=0}^\infty f_{k} \mu^{k(k-2)} z^{k}
     =\alpha \sum_{k=0}^\infty f_{k} \mu^{k(k-1)} (\mu^{-1} z)^{k}
     \equalscolon G(\overline\mu z).
 \end{multline*}
Consequently, the following is true.
\begin{corollary}\label{cr:ca2}\navy
    Given a real entire function~$\tilde F(z) = \sum_{k=0}^\infty f_kz^k$, $f_0>0$, of genus no
    more than one consider the functions
    \[
    G(z) = \alpha \sum_{k=0}^\infty i^{\frac{k(k-1)}2} f_{k} z^{k} \quad\text{and}\quad
    \overline{G}(z) = \overline{\alpha} \sum_{k=0}^\infty i^{-\frac{k(k-1)}2} f_{k} z^{k},
    \]
    where $\alpha$ is an arbitrary complex number not equal to zero.

    If the functions $\sum_{k=0}^{\infty} f_{2k}z^k$ and $\sum_{k=0}^{\infty} f_{2k+1}z^k$ have
    only negative zeros then for $\mu=\sqrt{i}$ the function $G(\overline\mu z)$ satisfies
    Theorem~\ref{th:ca1} with $c=\overline\mu$ and the function $\overline G(\mu z)$ satisfies
    Theorem~\ref{th:ca1} with $c=\mu$.
\end{corollary}

A real entire function $f(z^2) + z g(z^2)$ is \emph{strongly stable} if
$f(z^2) + (1+\zeta) z g(z^2)$ has no zeros in the closed right half of the complex plane for all
complex~$\zeta$ which are small enough. The notion of strong stability is not widely used. It is
the ``proper'' extension to entire functions of the polynomial stability. The strongly stable
functions include all stable polynomial and functions which are in some sense close to them (for
example, $e^x$); but they exclude such functions, like $e^{-x}$. For the case of complex
functions and further details see~\cite[p.~129]{Postnikov} or \emph{e.g.} the class~\textbf{HB}
in~\cite{ChebMei}.

Note that if $\tilde F(z)$ is strongly stable then the condition of Corollary~\ref{cr:ca2} is
satisfied. As it is noted in the introduction (and this is a consequence of Chebotarev's
Theorem, see~\cite{ChebMei,Postnikov}), strict stability of $\tilde F$ implies that the zeros of
$\sum_{k=0}^{\infty} f_{2k}z^k$ and $\sum_{k=0}^{\infty} f_{2k+1}z^k$ are negative, simple and
interlacing. At the same time, the interlacing property is redundant for this corollary.

\subsection{Meromorphic functions}

A careful observer will notice that Theorem~\ref{th:ca1} holds for any meromorphic function of
a certain form. Indeed, the following generalization is valid.
\begin{corollary}\label{cr:ca4}\navy
    Given an arbitrary entire function $E(z)$ suppose that a meromorphic function~$F(z)$ is such that
    $F(z)=f(z^2)+zg(z^2)$, where
    \[
    f(z) = f_0 z^{-k}e^{E(z)}\frac{\prod_{\mu=0}^{\omega_1}\left(1+\frac{z}{b_\mu}\right)}
                    {\prod_{m=0}^{\omega_3}\left(1-\frac{z}{\beta_\mu}\right)}
    \quad\text{and}\quad
    g(z) = g_0 z^je^{E(z)}\frac{\prod_{\nu=0}^{\omega_2}\left(1-\frac{z}{a_\nu}\right)}
                       {\prod_{n=0}^{\omega_4}\left(1+\frac{z}{\alpha_\nu}\right)},
    \]
    $f_0\ne 0$, $g_0\ne 0$, $j$ and $k$ are nonnegative integers,
    $b_\mu,a_\nu,\beta_m,\alpha_n>0$ for all $\mu=0,\dots,\omega_1$, $\nu = 0,\dots,\omega_2$,
    $m=0,\dots,\omega_3$, $n = 0,\dots,\omega_4$,
    $1\le\omega_1,\omega_2,\omega_3,\omega_4\le \infty$ and
    $\sum_{\mu=0}^{\omega_1} b_\mu^{-1} + \sum_{\nu=0}^{\omega_2} a_\nu^{-1} +
    \sum_{m=0}^{\omega_3}\beta_m^{-1} + \sum_{n=0}^{\omega_4}\alpha_n^{-1} < +\infty$.

    Then the zeros of $F(z)$ have the same localization as in Theorem~\ref{th:ca1}.
\end{corollary}

In this corollary the function $F(z)$ has the same zeros with multiplicities as the entire function
\[
    \tilde{F}(z) = f_0 \prod_{\mu=0}^{\omega_1}\left(1+\frac{z^2}{b_\mu}\right)
                       \prod_{n=0}^{\omega_4}\left(1+\frac{z^2}{\alpha_\nu}\right)
          + g_0 z^{2(j+k)+1} \prod_{m=0}^{\omega_3}\left(1-\frac{z^2}{\beta_\mu}\right)
                           \prod_{\nu=0}^{\omega_2}\left(1-\frac{z^2}{a_\nu}\right),
\]
which is satisfying Theorem~\ref{th:ca1}.

For the general case presented in Corollary~\ref{cr:ca4} the corresponding condition on the
coefficients of the power series is not so straightforward. This is because of arising
exponential factors in the functions $f(z)$ and~$g(z)$ defined by~\eqref{eq:f_even_odd}.

Furthermore, in fact the zeros of $F(z)$ must be the so-called
$c$-points\footnote{\emph{I.e.} the solutions to $R(z)=c$.} of the function $R(z)$; the double
(resp. triple) zeros of $F(z)$ are those $c$-points of $R(z)$ which are the points of
ramification\footnote{That is, the points satisfying $R'(z)=0$. The order of a ramification
    point is the multiplicity of this zero.}
of the first (resp. second) order. For the details see the formul\ae~\eqref{eq:cond_double_zero}
and~\eqref{eq:cond_triple_zero}.

\section{Proof of Theorem~\ref{th:ca1}}
\subsection{Basic inequalities}
We are interested in the points $z_*$ satisfying the equation~$F(z_*)=0$, which is equivalent
to~$R(z_*)=c$. Let us start calculating the logarithmic derivative of~$R(z)$.
\begin{multline}\label{DlogR}
    \frac{R'(z)}{R(z)} = (\ln R(z))'
    =\frac{\displaystyle\sum_{k=0}^{\omega_1}\frac{2z}{b_k}
        \prod_{\substack{\mu=0\\\mu\ne k}}^{\omega_1}
        \left(1+\frac{z^2}{b_\mu}\right)}
    {\displaystyle\prod_{\mu=0}^{\omega_1}
        \left(1+\frac{z^2}{b_\mu}\right)}
    - \frac {(2j+1)z^{2j}}{z^{2j+1}}
    -\frac{\displaystyle z^{2j+1}\sum_{k=0}^{\omega_2}\frac{2z}{-a_k}
        \prod_{\substack{\nu=0\\\nu\ne k}}^{\omega_2}
        \left(1-\frac{z^2}{a_\nu}\right)}
    {\displaystyle z^{2j+1}\prod_{\nu=0}^{\omega_2}\left(1-\frac{z^2}{a_\nu}\right)}
    \\
    = \sum_{k=0}^{\omega_1}\frac{2z}{b_k\left(1+\frac{z^2}{b_k}\right)}
    - \frac {2j+1}{z}
    + \sum_{k=0}^{\omega_2}\frac{2z}{a_k\left(1-\frac{z^2}{a_\nu}\right)}
    = \sum_{k=0}^{\omega_1}\frac{2z}{b_k+z^2}
    + \sum_{k=0}^{\omega_2}\frac{2z}{a_k-z^2}
    - \frac {2j+1}{z}.
\end{multline}

If we have $z=re^{i\phi}$, where $r>0$ and $-\pi<\phi\le\pi$, then
\[
\frac{\partial z}{\partial r} = e^{i\phi} = \frac{z}{|z|}
\quad\text{and}\quad 
\frac{\partial z}{\partial\phi} = ire^{i\phi} = iz.
\]
Therefore,
\begin{equation}\label{eq:arg_R}
    \frac{\partial}{\partial r} \arg R(z) = \frac {\partial}{\partial r} \Im\ln R(z)
    = \Im \frac {\partial\ln R(z)}{\partial r} = \Im \left(\frac{z}{|z|} (\ln R(z))'\right)
    = \frac 1{|z|} \Im \big(z(\ln R(z))'\big),
\end{equation}
and
\begin{equation}\label{eq:ln_abs_R}
    \frac{\partial}{\partial\phi} \ln|R(z)| = \frac {\partial}{\partial\phi} \Re\ln R(z)
    = \Re \frac {\partial\ln R(z)}{\partial\phi} = \Re \left(iz(\ln R(z))'\right)
    = -\Im \big(z(\ln R(z))'\big),
\end{equation}
The coincidence on the right-hand side of these two equalities reflects the Cauchy-Riemann
condition. Observe, for all passible $k$
\begin{align*}
    \Im \frac{2z^2}{a_k-z^2}
    = \Im \frac{2}{a_kz^{-2}-1} &> 0&&\text{and}&
    \Im \frac{2z^2}{b_k+z^2}
    = \Im \frac{2}{b_kz^{-2}+1} &> 0&&\text{whenever}&\Im z^2&>0,\\
    \Im \frac{2z^2}{a_k-z^2} &< 0&&\text{and}& 
    \Im \frac{2z^2}{b_k+z^2} &< 0&&\text{whenever}&\Im z^2&<0.
\end{align*}
The summation of these inequalities on $k$, in accordance with~\eqref{DlogR}, gives
\begin{equation}\label{eq:dR_sign}
    \sign \Im \big(z(\ln R(z))'\big)
    = \sign\Im\left(\sum_{k=0}^{\omega_1}\frac{2z^2}{b_k+z^2}
        + \sum_{k=0}^{\omega_2}\frac{2z^2}{a_k-z^2}
        - \frac {(2j+1)z}{z}\right)
    = \sign\Im z^2.
\end{equation}
Consequently, the identities~\eqref{eq:arg_R} and~\eqref{eq:ln_abs_R} imply
\begin{align}\label{eq:dR_monotonic_r}
    &\sign\frac{\partial}{\partial r} \arg R(z)
    =\sign\Im z^2,\quad\text{and}
    \\\label{eq:dR_monotonic_phi}
    &\sign\frac{\partial}{\partial\phi} \ln|R(z)|
    =-\sign\Im z^2.
\end{align}

\subsection{The geometry of the set $\ln|R(z)|=0$}
\label{sec:the-geom-set}
For any harmonic function $u(z)$ which is not a constant the following is true (see,
\emph{e.g.}~\cite[pp.~18--19]{Duren}):
\begin{enumerate}
\item \label{item:harmonic_prop1}
    it can only have isolated critical points, i.e. points where
    \[
    \frac{\partial u}{\partial z} \colonequals \frac 12\left(\frac{\partial u}{\partial x} - i
        \frac{\partial u}{\partial y}\right) = 0,\quad\text{where}\ x=\Re z,\ y=\Im z;
    \]
\item \label{item:harmonic_prop2}
    the level set $u(z)=u(z_0)$ in a neighbourhood of any point~$z_0$ not critical for $u(z)$ is an
    analytic arc;
\item \label{item:harmonic_prop3}
    the level set $u(z)=u(z_0)$ in a neighbourhood of a critical point~$z_0$ is the
    intersection of $n$ analytic arcs under the angles $\frac{\pi}n$. Here $n-1$ is the
    multiplicity of the zero at~$z_0$ of the analytic function~$\frac{\partial u}{\partial z}$.
\end{enumerate}

Denote
\[
u(z)\colonequals\ln|R(z)|,\quad
\Gamma\colonequals\{z\in\mathbb{C}:u(z)=0\}
\quad\text{and}\quad \Gamma_i\colonequals\Gamma\cap Q_i,\ i=1,2,3,4.
\]
The function $u(z)$ is harmonic on $\mathbb{C}\setminus(A\cup B)$,
approaching~$A\colonequals\{0,\pm\sqrt{a_1},\pm\sqrt{a_2},\dots\}$ it tends to $+\infty$, and
approaching~$B\colonequals\{\pm i\sqrt{b_1},\pm i\sqrt{b_2},\dots\}$ it tends to $-\infty$. Therefore,
$u(z)\ne 0$ close to singulariries, i.e. there are no singulariries of~$u(z)$ close to~$\Gamma$.

From~\eqref{eq:dR_monotonic_phi} it follows that $\frac{\partial u(z)}{\partial z}\ne 0$ for any
$z$ satysfying~$\Im z^2\ne0$. That is, all critical points of $u(z)$ (if any) are on the complex
and real axes (see~\ref{item:harmonic_prop1}).


Observe that
\begin{equation}\label{eq:R_conj_neg}
    R(z)
    =\frac{1}{z^{2j+1}}\cdot
            \frac{|f_0|\prod_{\mu=0}^{\omega_1}\left(1+\frac{z^2}{b_\mu}\right)}
                 {|g_0|\prod_{\nu=0}^{\omega_2}\left(1-\frac{z^2}{a_\nu}\right)}
    =\overline{R(\overline{z})}=-R(-z),
\end{equation}
which implies
\begin{equation}\label{eq:u_symm}
    u(z)=u(\overline{z})=u(-z)=u(-\overline{z}).
\end{equation}
This condition, in particular, provides that for $i=1,2,3$ the set~$\Gamma_{i+1}$ can be
obtained from~$\Gamma_{i}$ by reflecting from the semi-axis
$\overline{Q}_i\cap\overline{Q}_{i+1}$.

Suppose that $z\in Q_1$, then $\Im z^2> 0$ and $u(z)$ is a monotonically decreasing function of
$\phi=\Arg z$ by the condition~\eqref{eq:dR_monotonic_phi}. Therefore, for a fixed $r>0$ the
function~$u(re^{i\phi})$ reaches~$0$ in~$\overline{Q}_1$ at most for one value
$\phi\in[0,\frac\pi2]$. In particular, this implies that the value of~$n$ (the number of
intersecting arcs) in the property~\ref{item:harmonic_prop3} is at most~$2$. Indeed, consider a
critical point $z_0$ on $\overline{Q}_1\setminus Q_1$; without loss of generality assume that
$z_0\in\overline{Q}_1\cap\overline{Q}_4$. If $n=4$ or more, then in a sufficiently small
neighbourhood of $z_0$ there must be at least two points with the same absolute value satisfying
$u(z)=0$, which is not the case. If $n=3$ then at least in one of the quadrants $Q_1$, $Q_4$
there must be at least two points with the same absolute value satisfying $u(z)=0$, which by the
condition~\eqref{eq:u_symm} will imply the same for the remaining quadrant.

We conclude
\emph{\navy that\/ $\Gamma\cap\overline{Q}_1 = \overline{\Gamma}_1$, and this set is a union of
    closed analytic arcs. These arcs can meet only on\/~$\overline{\Gamma}_1\setminus\Gamma_1$,
    and there are at most two arcs which can meet at one point. Different points
    of\/~$\overline{\Gamma}_1$ have different absolute values, and this is also valid at
    critical points (if any).}
These facts are true for $\overline{\Gamma}_2$, $\overline{\Gamma}_3$ and $\overline{\Gamma}_4$
as well, provided by the condition~\eqref{eq:u_symm}.

\subsection{The value of $\arg R(z)$ on $\overline{\Gamma}_1$}
\label{sec:arg-R}
Since from~\eqref{eq:dR_monotonic_phi} we have that $u(re^{i\phi})$ is decreasing in
$\phi\in(0,\frac \pi 2)$ for a fixed $r$,  
\begin{equation}\label{eq:ln_in_between}
    u(ir) = u(re^{i\frac{\pi}{2}})<u(re^{i\phi})<u(re^0) = u(r),
    \quad\text{for all }\phi\in(0,\tfrac{\pi}2).
\end{equation}
Taking into account the conclusion (of Subsection~\ref{sec:the-geom-set}) from the previous
paragraph, we note that the function $l:(0,+\infty)\to[0,\frac{\pi}2]$ such that
\[
l(r)\colonequals\begin{cases}
    0,&\text{if }u(r)\le 0,\\
    \frac{\pi}2,&\text{if }u(ir)\ge 0,\\
    \phi&\text{satisfying }u(re^{i\phi})=0 \text{ otherwise,}
\end{cases}
\]
is defined correctly. Moreover, this function is piecewise-analytic and continuous. Indeed,
consider any $s>0$. Then if there is a $\phi\in[0,\frac\pi2]$ solving the equation
$u(se^{i\phi})=0$ (such $\phi$ can be only unique), then $l(s)=\phi$ by definition. In this
case, $l(r)$ is a parametrization of some analytic arc of $\Gamma_1$. If $u(se^{i\phi})\ne 0$
for all $\phi\in[0,\frac\pi2]$ then $l(s)=0$ or $l(s)=\frac\pi2$. If so, $u(se^{il(s)})$ being
not singular has the same sign for any $r$ close enough to $s$, and hence $l(r)=l(s)$. That is,
in the neighbourhood of $s$ the function $l(r)$ is constant. Noting that such intervals can be
bounded only by the points where~$u(re^{il(r)})$ vanishes (and by the origin), we deduce that
$l(r)$ is continuous for all $r>0$.

Considering the set $\{re^{il(r)}:r>0\}\setminus\Gamma_1$ allows us to deduce that $\arg R(z)$
has no ``jumps'' when we pass from one arc of $\Gamma_1$ to another. We have that $\arg R(z)$
and $u(z)$ are non-singular harmonic functions at the same points $z\in\mathbb{C}$. So, in
particular, \emph{\navy $\arg R(re^{il(r)})$ depends continuously on $r$}; i.e. we can take a
branch of $\arg R(re^{il(r)})$ which is continuous in $r$.

The condition~\eqref{eq:dR_monotonic_r} implies that $\arg R(re^{i\phi})$ for $r>0$ increases
(strictly) for a fixed $\phi\in(0,\frac\pi2)$ and remains constant if $\phi=0$ or
$\phi=\frac\pi2$. Given $z\in\Gamma_1$ let us consider the local coordinates~$(\tau,\nu)$ about
the point $P\colonequals z$ such that the axis~$P\tau$ is tangent to~$\Gamma_1$ at~$P$, and the
axis~$P\nu$ is orthogonal to $\Gamma_1$ at $P$.

We orient the curve~$\Gamma_1$ such that its positive direction corresponds to the growth of
$|z|$ (\emph{i.e.} towards the infinity from the origin~$z=0$). This is possible (see
Subsection~\ref{sec:the-geom-set}) and corresponds to the condition
\[
\frac{\partial\tau}{\partial r}>0.
\]
The chain rule and the Cauchy-Riemann condition give us
\begin{equation*}
    0 < \frac{\partial \arg R(z)}{\partial r}
    = \frac{\partial \arg R(z)}{\partial\tau}\frac{\partial\tau}{\partial r}
    + \frac{\partial \arg R(z)}{\partial\nu}\frac{\partial\nu}{\partial r}
    = \frac{\partial \arg R(z)}{\partial\tau}\frac{\partial\tau}{\partial r}
    \pm \frac{\partial u(z)}{\partial\tau}\frac{\partial\nu}{\partial r}
    = \frac{\partial \arg R(z)}{\partial\tau}\frac{\partial\tau}{\partial r}
\end{equation*}
Therefore, $\frac{\partial \arg R(z)}{\partial\tau}>0$ everywhere on~$\Gamma_1$. So we showed
that
\emph{\navy $\arg R(re^{il(r)})$ is a continuous non-decreasing function of~$r>0$; it increases if and
    only if\/ $0<l(r)<\frac\pi2$.}

Since about the origin $R(z)\sim\frac{|f_0|}{|g_0|z^{2j+1}}$, $j\ge0$,
\emph{\navy for all~$r>0$ small enough we have $l(r)=\frac\pi2$ and}
\begin{equation}\label{eq:arg_at_O}
    {\navy \Arg R(re^{il(r)})=\Arg R(ir)= \Arg i^{-2j-1} = (-1)^{j+1}\frac\pi2.}
\end{equation}
Furthermore, according to the condition~\eqref{eq:dR_monotonic_phi},
\begin{equation}\label{eq:min_r_arg}
\textit{\navy for the minimal $r>0$ satisfying $u(re^{il(r)})=0$ we have $l(r)=(-1)^{j+1}\dfrac\pi2$.}
\end{equation}

\subsection{$F(z)$ can only have simple zeros outside of the real and imaginary axes}
\label{sec:F-simple-zeros}
Note that if $R=p/q$, then
\begin{equation}\label{eq:dR}
    R'= \left(\frac {p}{q}\right)'
    = \left(\frac {p'}{q} - \frac {pq'}{q^2} \right)
    = \frac{1}{q}\left(p' - \frac {p}{q}q'\right)
    = \frac{p' - R q'}{q},
\end{equation}
\begin{multline}\label{eq:ddR}
    R'' = \left(\frac{p' - R q'}{q}\right)'
    = \frac {p''q - R'qq' - R qq'' - q'p' + Rq'^2}{q^2}\\
    = \frac {p'' - R q''}{q}- R'\frac {q'}{q} - \frac{q'(p'-Rq')}{q^2} = \frac {p'' - R q''}{q}
    - 2R'\frac {q'}{q}.
\end{multline}

The equality~\eqref{eq:dR_sign} shows that~$F(z)$ has no multiple zeros on
$\{z\in\mathbb{C}:\Im z^2\ne 0\}$. Indeed, let $F(z_0)=f(z_0^2)+z_0g(z_0^2)=0$ and
$\Im z_0^2\ne 0$. Then $R(z_0) = c$, so from~\eqref{eq:dR} we have
\begin{equation*}
    R'(z_0)
    = - \frac{(-cf(z_0^2))' - c (zg(z_0^2))'}{z_0 g(z_0^2)}
    = - c\frac{F'(z_0)}{z_0 g(z_0^2)}.
\end{equation*}
Consequently, if $F(z)$ has a multiple zero at~$z=z_0$, i.e. $F(z_0)=F'(z_0)=0$, then
\begin{equation}\label{eq:cond_double_zero}
\frac{R'(z_0)}{R(z_0)} = \frac 1c R'(z_0) = - \frac{F'(z_0)}{z_0 g(z_0^2)} = 0,
\end{equation}
which is contradicting to~\eqref{eq:dR_sign}.

\subsection{$F(z)$ can only have double and simple zeros on the real and imaginary axes}
Suppose that we have a multiple (at least triple) root at the point~$z_0\in\mathbb{C}$. Then
\begin{equation*}
    R(z_0) = c,\quad R'(z_0)=0\quad\text{from \eqref{eq:cond_double_zero},}\quad\text{and}\quad F''(z_0)=0.
\end{equation*}
Therefore, substituting of $p(x_0)=f(x_0^2)$ and $q(x_0)=x_0g(x_0^2)$ into~\eqref{eq:ddR} gives
us
\begin{equation}\label{eq:cond_triple_zero}
    R''(x_0)
    = \frac {-c(f(x_0^2))'' - c(x_0g(x_0^2))''}{x_0g(x_0^2)}
    = \frac {-cF''(x_0)}{x_0g(x_0^2)} = 0
\end{equation}
As a consequence,
\[
(\ln R(z))''|_{z=x_0} =\left.\left(\frac{R'(z)}{R(z)}\right)'\right|_{z=x_0}
=\frac {R''(x_0)}{R(x_0)} - \frac {(R'(x_0))^2}{R^2(x_0)}=0.
\]

Then the harmonic function $u(z)=\ln|R(z)|$ satisfies
\[
\frac{\partial u(z)}{\partial z}= (\ln R(z))'=0,\quad
\frac{\partial^2 u(z)}{\partial z^2}= (\ln R(z))''=0.
\]
However, this contradicts the conclusion of the subsection~\ref{sec:the-geom-set}. Therefore,
\emph{\navy the multiplicity of any zero of $F$ does not exceed $2$.} Note, that each double
zero must satisfy~$\Im z_0^2=0$ (see Subsection~\ref{sec:F-simple-zeros}).

\subsection{Distribution of zeros of $F(z)$}
All zeros of~$F(z)$ from the first quadrant~$\overline Q_1$ are in the set~$\overline\Gamma_1$.
The set~$\overline\Gamma_1$ in its turn is a subset of $\{re^{il(r)}:r>0\}$ (see
Subsection~\ref{sec:arg-R}). Therefore, increasing positive~$r$ from~$0$ to infinity we obtain
the series of all solutions from $\overline{Q}_1$ to the following problems:
\begin{align}
    &\ &&\ln|R(z)|=1,&& \arg R(z)=&\arg c&&
    \xLongleftrightarrow{\qquad} &&R(z)=&&c,&&\ &&\label{eq:probl_Q1}\\
    &\ &&\ln|R(z)|=1,&& \arg R(z)=&\pi-\arg c&&
    \xLongleftrightarrow{\qquad} &&R(z)=&&-\overline{c}.&&\ &&\label{eq:probl_Q2}\\
    &\ &&\ln|R(z)|=1,&& \arg R(z)=&\pi+\arg c&&
    \xLongleftrightarrow{\qquad} &&R(z)=&&-c,&&\ &&\label{eq:probl_Q3}\\
    &\ &&\ln|R(z)|=1,&& \arg R(z)=&-\arg c&&
    \xLongleftrightarrow{\qquad} &&R(z)=&&\overline{c},&&\ &&\label{eq:probl_Q4}
\end{align}
On the other hand, according to~\eqref{eq:R_conj_neg},
\begin{align*}
    \arg R(\overline{z}) &= \arg \overline{R}(z) = -\arg R(z),\\
    \arg R(-z) &= \arg (-R(-z)) = \pi+\arg R(z),\\
    \arg R(-\overline{z}) &= \arg (-\overline{R}(z)) = \pi-\arg R(z)
\end{align*}
for all $z\in Q_1$. As a consequence,
\begin{itemize}
\item the equality~\eqref{eq:probl_Q1} gives all zeros of $F(z)$ in $\overline Q_1$,
\item the equality~\eqref{eq:probl_Q2} gives all zeros of $F(z)$ in $\overline Q_2$,
\item the equality~\eqref{eq:probl_Q3} gives all zeros of $F(z)$ in $\overline Q_3$,
\item the equality~\eqref{eq:probl_Q4} gives all zeros of $F(z)$ in $\overline Q_4$.
\end{itemize}
\begin{figure}[h]
    \centering
    \setlength{\unitlength}{1mm}
    \begin{picture}(194,62)(-24,1)
        \thinlines
        \put(0,10){\vector(1,0){60}}
        \put(55,5){\text{$\Re$}}
        \put(10,0){\vector(0,1){60}}
        \put(12,55){\text{$\Im$}}
        \put(6,5){\text{$0$}}
        \green
        \thicklines
        \put(10,10){\line(0,1){5}}
        \qbezier(10,15)(19,15)(22,10)
        \qbezier(29,10)(29,11)(29.4,13.5)
        \qbezier(29.4,13.5)(30.5,17)(30.2,20)
        \qbezier(30.2,20)(29.5,27)(16,33)
        \qbezier(16,33)(13.6,34.1)(12,36.5)
        \qbezier(12,36.5)(11,38)(10,39)
        \qbezier(10,39)(12,41)(14,44)
        \qbezier(14,44)(16,47)(19,45.5)
        \qbezier(19,45.5)(27,42)(40,40)
        \qbezier(40,40)(49,39)(55,48)
        \qbezier(55,48)(59,54.5)(60,60)
        \put(17.5,5.8){\text{$u(z)=0$}}
        \gray
        \thinlines
        
        \qbezier(10,10)(16.0,11.5)(17,14)
        \qbezier(17,14)(19.5,20)(12,24)
        \qbezier(12,24)(10.5,25.1)(10,26.4)
        \put(-9.6,20.35){\text{$\frac{R(z)}{|R(z)|}=\overline{c}$}}
        \put(7,21.15){\vector(3,1){6}}
        \qbezier(42,10)(37.5,11)(35.5,17)
        \qbezier(35.5,17)(33.5,26.2)(28.7,22.8)
        \qbezier(28.7,22.8)(24,20.5)(20.5,24)
        \qbezier(20.5,24)(13.5,32)(10,26.4)
        \put(-9.6,31.0){\text{$\frac{R(z)}{|R(z)|}={c}$}}
        \put(7,31.3){\vector(3,-2){4.5}}
        \qbezier(42,10)(38,14)(36,25)
        \qbezier(36,25)(34,34)(23,37)
        \qbezier(23,37)(14,40)(10,45)
        \put(-12.85,40.65){\text{$\frac{R(z)}{|R(z)|}=-\overline{c}$}}
        \put(7,41.5){\vector(3,1){4.5}}
        %
        \qbezier(42,10)(57,30)(45,45.5)
        \qbezier(45,45.5)(37.8,54.3)(21,51.5)
        \qbezier(21,51.5)(13,50)(10,45)
        \put(-12.85,49.8){\text{$\frac{R(z)}{|R(z)|}=-{c}$}}
        \put(7,50.1){\vector(2,-1){4.8}}
        \qbezier(42,10)(58,16)(59,36)
        \qbezier(59,36)(59.6,47)(48,52)
        \qbezier(48,52)(39,56)(35,60)
        \put(35.95,3.8){\text{$\frac{R(z)}{|R(z)|}=\overline{c}$}}
        \put(44.4,5.75){\vector(1,2){3.5}}
        \bla
        \put(10,10){\circle*{.8}}
        \put(13.5,16){\text{$\overline{z}_1$}}
        \put(29.2,24.75){\text{$z_2$}}
        \put(14.1,41.8){\text{$-\overline{z}_3$}}
        \put(42,37.2){\text{$-z_4$}}
        \put(56.6,47){\text{$\overline{z}_5$}}
        \red
        \put(16.95,13.9){\circle*{.8}}
        \put(29.2,23.05){\circle*{.8}}
        \put(12.8,42.2){\circle*{.8}}
        \put(47.7,41.2){\circle*{.8}}
        \put(54.75,47.6){\circle*{.8}}
        %
        %
        \bla
        \thinlines
        \put(80,10){\vector(1,0){60}}
        \put(135,5){\text{$\Re$}}
        \put(90,0){\vector(0,1){60}}
        \put(92,55){\text{$\Im$}}
        \put(86,5){\text{$0$}}
        \green
        \thicklines
        \put(90,10){\line(0,1){6}}
        \put(90,16){\line(1,0){1}}
        \qbezier(91,16)(96.2,17)(97.2,20.9)
        \qbezier(97.2,20.9)(98.5,25.1)(104,27.6)
        \qbezier(104,27.6)(110.0,29.9)(114.8,37.4)
        \qbezier(114.8,37.4)(115.9,39.0)(117.45,42.0)
        \qbezier(117.45,42.0)(119.0,45.0)(121.0,47.35)
        \qbezier(121.0,47.35)(126.5,53.0)(131.2,55)
        \qbezier(131.2,55)(132.5,55.8)(135,56.8)
        \qbezier(135,56.8)(139.15,58.9)(140,59.9)
        \put(104,46.2){\text{$u(z)=0$}}
        \gray
        \thinlines
        \qbezier(90,10)(95.9,15.9)(92.5,18)
        \qbezier(92.5,18)(90.9,18.7)(90,20)
        \put(70.4,15.4){\text{$\frac{R(z)}{|R(z)|}=\overline{c}$}}
        \put(87,16.3){\vector(3,1){5}}
        \qbezier(90,20)(99,24.3)(112,12)
        \qbezier(112,12)(113,11)(114.5,10)
        \put(70.4,25.9){\text{$\frac{R(z)}{|R(z)|}={c}$}}
        \put(87,26.0){\vector(1,-1){5}}
        \qbezier(114.5,10)(115.6,23.0)(100,37.9)
        \qbezier(100,37.9)(96.7,41.3)(90,45.5)
        \put(67.15,38.3){\text{$\frac{R(z)}{|R(z)|}=-\overline{c}$}}
        \put(87,39.4){\vector(2,1){6.7}}
        \qbezier(90,45.5)(104.2,47.0)(118.0,36)
        \qbezier(118.0,36)(122.3,32.8)(128.5,26)
        \qbezier(128.5,26)(137.2,15.5)(138,10)
        \put(67.15,53.1){\text{$\frac{R(z)}{|R(z)|}=-{c}$}}
        \put(87,53.1){\vector(2,-3){4.8}}
        \qbezier(138,10)(143,24)(136,40)
        \qbezier(136,40)(130.7,51)(124,60)
        \put(106.4,54.5){\text{$\frac{R(z)}{|R(z)|}=\overline{c}$}}
        \put(123.1,55.2){\vector(1,0){4}}
        %
        \bla
        \put(90,10){\circle*{.8}}
        \put(95.1,13.8){\text{$\overline{z}_1$}}
        \put(93.3,23.0){\text{$z_2$}}
        \put(102.7,24.8){\text{$-\overline{z}_3$}}
        \put(116.1,37.2){\text{$-z_4$}}
        \put(130.6,51.3){\text{$\overline{z}_5$}}
        \red
        \put(93.55,16.75){\circle*{.8}}
        \put(97.2,20.9){\circle*{.8}}
        \put(107.55,29.45){\circle*{.8}}
        \put(115.22,38){\circle*{.8}}
        \put(128.45,53.5){\circle*{.8}}
    \end{picture}
    \caption{Examples of zero location on the curve~$u(z)=0$ for a function~$F(z)$ in the case
        when $c\in Q_1$ and $j=1$ (a scheme). See also a calculated plot of a sample polynomial in
        the appendix on the page~\pageref{fig:3}.}
\label{fig:0}
\end{figure}
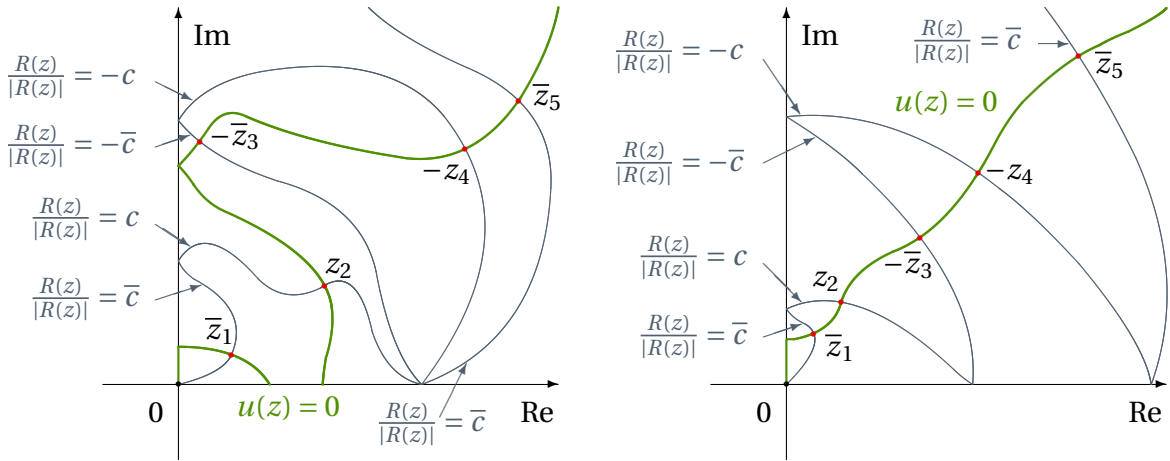

Now we only have to arrange these solutions. Note that for $\Im c^2\ne 0$ the function
$\arg R(z)$ never achieves the value $\arg c$ when $\Im z^2\ne 0$.
\begin{SCfigure}[1.15][h]
    \centering
    \setlength{\unitlength}{1mm}
    \begin{picture}(72,62)(-9,-1)
        \thinlines
        \put(0,30){\vector(1,0){60}}
        \put(55,25){\text{$\Re$}}
        \put(30,0){\vector(0,1){60}}
        \put(32,55){\text{$\Im$}}
        \put(26,25){\text{$0$}}
        \green
        \thicklines
        \put(30,30){\line(0,-1){5}}
        \qbezier(30,25)(35,25)(38,30)
        \qbezier(38,30)(41,38)(34,42)
        \qbezier(34,42)(28.7,44.8)(23,42)
        \qbezier(23,42)(16,38)(14,30)
        \put(14,30){\line(-1,0){4}}
        \qbezier(10,30)(10,20)(14,15)
        \qbezier(14,15)(19,8)(30,5)
        \qbezier(30,5)(47,2.5)(56,16)
        \put(40,2){\text{$z=r\frac{R(re^{l(r)})}{|R(re^{l(r)})|}$}}
        \gray
        \thinlines
        \put(30,30){\line(3,2){27}}
        \put(43.5,48){\text{$\frac {z}{|z|}=c$}}
        \put(30,30){\line(-3,-2){27}}
        \put(2,8.5){\text{$\frac {z}{|z|}=-c$}}
        \put(30,30){\line(3,-2){27}}
        \put(36.2,15){\text{$\frac {z}{|z|}=\overline c$}}
        \put(30,30){\line(-3,2){27}}
        \put(6.5,48){\text{$\frac {z}{|z|}=-\overline c$}}
        \bla
        \put(30,30){\circle*{.8}}
        \put(34,23.6){\text{$z_1$}}
        \put(39.5,34.3){\text{$z_2$}}
        \put(15.7,40.1){\text{$z_3$}}
        \put(7.7,18.3){\text{$z_4$}}
        \put(53.2,10.5){\text{$z_5$}}
        \red
        \put(35.1,26.6){\circle*{.8}}
        \put(38.8,35.8){\circle*{.8}}
        \put(18.05,37.95){\circle*{.8}}
        \put(12.1,18.05){\circle*{.8}}
        \put(54.3,13.65){\circle*{.8}}
    \end{picture}
    \caption{An example of zero location on the curve\newline
        \(\displaystyle z=re^{i\arg R(re^{l(r)})}=r\frac{R(re^{l(r)})}{|R(re^{l(r)})|}\)\newline
        for a function~$F(z)$ in the case when $c\in Q_1$ and $j$ is even (a scheme).}
\label{fig:1}
\end{SCfigure}
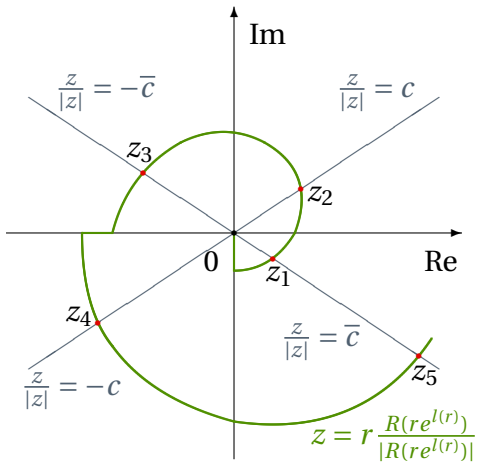
Suppose in the expressions~\eqref{eq:R_form} and~\eqref{eq:arg_at_O} the integer~$j$ is even, then
for small~$r>0$ we have $\Arg R(re^{il(r)}) = -\frac\pi2$. For $\Im c^2\ne 0$ the zeros
$\{z_k\}_{k=1}^\infty$, $0<|z_1|<|z_2|<|z_3|<\dots$, of $F(z)$ are simple and distributed as
follows:
\begin{itemize}
\item if $c\in Q_1$, then $z_1\in Q_4$, $z_2\in Q_1$, $z_3\in Q_2$, $z_4\in Q_3$, $z_5\in Q_4$
    and so on;
\item if $c\in Q_2$, then $z_1\in Q_3$, $z_2\in Q_2$, $z_3\in Q_1$, $z_4\in Q_4$, $z_5\in Q_3$
    and so on;
\item if $c\in Q_3$, then $z_1\in Q_2$, $z_2\in Q_3$, $z_3\in Q_4$, $z_4\in Q_1$, $z_5\in Q_2$
    and so on;
\item if $c\in Q_4$, then $z_1\in Q_1$, $z_2\in Q_4$, $z_3\in Q_3$, $z_4\in Q_2$, $z_5\in Q_1$
    and so on;
\end{itemize}

Given $k,l \in \{1,2,3,4\}$ denote $Q_{kl}\colonequals\overline Q_k \cap \overline Q_l$.  For
the remaining values of $c$ we use the condition~\eqref{eq:min_r_arg}, so the zeros
$\{z_k\}_{k=1}^\infty$ of $F(z)$, satisfy the following.
\begin{itemize}
\item if $c=1$, then $0<|z_1|\le|z_2|<|z_3|\le|z_4|<|z_5|\dots$ counting with multiplicities,\\
    $z_1\in Q_4\cup Q_{14}$ and $z_2\in Q_1\cup Q_{14}$ such that $\arg z_1=-\arg z_2$,\\
    $z_3\in Q_2\cup Q_{23}$ and $z_4\in Q_3\cup Q_{23}$ such that $\arg z_3=-\arg z_4$,\\
    $z_5\in Q_4\cup Q_{14}$ and $z_6\in Q_1\cup Q_{14}$ such that $\arg z_5=-\arg z_6$, and so
    on;
\item if $c=i$, then $0<|z_1|<|z_2|\le|z_3|<|z_4|\le|z_5|\dots$ counting with multiplicities,\\
    $z_1\in Q_{34}$,\\
    $z_2\in Q_2\cup Q_{12}$ and $z_3\in Q_1\cup Q_{12}$ such that $\arg z_2=\pi-\arg z_3$,\\
    $z_4\in Q_4\cup Q_{34}$ and $z_5\in Q_3\cup Q_{34}$ such that $\arg z_4=\pi-\arg z_5$,\\
    $z_6\in Q_2\cup Q_{12}$ and $z_7\in Q_1\cup Q_{12}$ such that $\arg z_6=\pi-\arg z_7$,\\
    and so on;
\item if $c=-1$, then $0<|z_1|\le|z_2|<|z_3|\le|z_4|<|z_5|\dots$ counting with multiplicities,\\
    $z_1\in Q_2\cup Q_{23}$ and $z_2\in Q_3\cup Q_{23}$ such that $\arg z_1=-\arg z_2$,\\
    $z_3\in Q_4\cup Q_{14}$ and $z_4\in Q_1\cup Q_{14}$ such that $\arg z_3=-\arg z_4$,\\
    $z_5\in Q_2\cup Q_{23}$ and $z_6\in Q_3\cup Q_{23}$ such that $\arg z_5=-\arg z_6$, and so
    on;
\item if $c=-i$, then $0<|z_1|<|z_2|\le|z_3|<|z_4|\le|z_5|\dots$ counting with multiplicities,\\
    $z_1\in Q_{12}$,\\
    $z_2\in Q_4\cup Q_{34}$ and $z_3\in Q_3\cup Q_{34}$ such that $\arg z_2=\pi-\arg z_3$,\\
    $z_4\in Q_2\cup Q_{12}$ and $z_5\in Q_1\cup Q_{12}$ such that $\arg z_4=\pi-\arg z_5$,\\
    $z_6\in Q_4\cup Q_{34}$ and $z_7\in Q_3\cup Q_{34}$ such that $\arg z_6=\pi-\arg z_7$,\\
    and so on;
\end{itemize}

Suppose in the expressions~\eqref{eq:R_form} and~\eqref{eq:arg_at_O} the integer~$j$ is odd,
then for small~$r>0$ we have $\Arg R(re^{il(r)}) = \frac\pi2$. Then the zeros of $F(z)$ are
distributed as in the case of even~$j$ after the substitution~$c\mapsto -c$.  The theorem is
proved completely.

\section{A conjecture by Alan Sokal}\label{sec:conj-alan-sokal}
\subsection{Statement of the conjecture}
\paragraph{ }
Alan Sokal in his Talk at Institut Henri Poincar\'e (9th of November 2009, \cite{Sokal}) put
forward the hypothesis that
\begin{conjecture}\label{conjecture-A}
    The entire function
    \begin{equation}\label{eq:0}
        F(z;q)=F(z)=\sum_{k=0}^\infty \frac {q^{\frac{k(k-1)}2}z^k}{k!},
    \end{equation} 
    where $q$ is a complex number, $0<|q|\leq 1$, can have only simple zeros.
\end{conjecture}
The stronger version of the conjecture claims that
\begin{conjecture}\label{conjecture-B}
    The function $F(z;q)$ for $q\in\mathbb{C}$, $\ 0<|q|\leq 1$, can have only simple zeros with
    distinct absolute values.
\end{conjecture}

\paragraph{ }
The following facts on~$F(z)$ are known. The function $F(z)$ is the unique solution to the
following Cauchy problem
\begin{equation}\label{eq:1}
F'(z)=F(qz),\quad F(0)=1,
\end{equation}
which can be checked directly. Moreover, when $|q|=1$ this function has the exponential
type~$1$, for~$q$ lower in modulus the function~$F(z)$ is of zero genus.

The case of positive~$q$ was studied extensively. It is shown that the zeros of $F$ are negative
(see~\cite{Morris_et_al}), simple and satisfy Conjecture~\ref{conjecture-B} as well as certain
further conditions (see~\cite{Liu,Langley}). In Appendix~\ref{sec:conjecture-real-q}, using
almost the same argumentation as in~\cite{Liu} we show that Conjecture~\ref{conjecture-B} holds
true for all real~$q$. However, we also obtain the same result using another technique (see
Subsection~\ref{sec:conj-negat-q}).

The properties of $F(z)$ for complex~$q$ was studied
in~\cite{Alander,Valiron,Eremenko_Ostrovsky}. According to~\cite{Sokal},
Conjecture~\ref{conjecture-B} is true if $|q|<1$ and the zeros of~$F(z)$ big enough in absolute
value (A.~Eremenko) as well as for small $|q|$.

\paragraph{ }
The family of polynomials
\begin{equation}\label{eq:01}
    P_N(z;q)=\sum_{k=0}^N \binom{N}{n}z^nq^{\frac{n(n-1)}{2}}
\end{equation} 
is connected tightly to the function $F(z)$; it approximates this function in the sense that
\[
P_N\left(zN^{-1};q\right)
    = \sum_{k=0}^N \frac{q^{\frac{n(n-1)}{2}}z^n}{n!}
           \left(1-\frac{n}N+\frac{n-1}N\right)\cdot\left(1-\frac{n}N+\frac{n-2}N\right)
               \cdots\left(1-\frac{n}N+\frac{1}N\right) \xrightarrow{N\to\infty} F(z;q).
\]
The polynomial version of this conjecture have the following form.
\begin{conjecture}\label{conjecture-C}
    For all $N>0$ the polynomial $P_N(z;q)$ where $|q|<1$, can have only simple roots, separated
    in modulus by at least the factor $|q|$.
\end{conjecture}
The original (equivalent) statement is concerned with the family of polynomials
$\left\{P_N\left(zw^{N-1},w^{-2}\right)\right\}_{N\in\mathbb{N}}$, where $w^{-2}=q$.

\paragraph{ }
Observe that $\ref{conjecture-C}\implies \ref{conjecture-B}\implies \ref{conjecture-A}$.  In the
present paper (see Subsections~\ref{sec:conj-negat-q} and \ref{sec:case-q2=-rho2}) we prove that
Conjecture~\ref{conjecture-B} is true whenever~$q^2\in\mathbb{R}$. The approach is translated
without changes on polynomials $P_N(z;q)$. This gives the assertion of
Conjecture~\ref{conjecture-C} without bounds on the ratio of subsequent (by modulus) roots.

\subsection{The properties of the even and the odd
    parts of~$F(z;\rho)$ for the positive~$\rho$}
\label{sec:properties-even-odd}

For a function~$\Psi(z)$ we define its even and odd parts as follows
\begin{equation}\label{eq:odd_even_def}
\Psi_e(z) \overset{\rm def}= \frac{\Psi(z)+\Psi(-z)}{2}
\qquad\text{and}\qquad
\Psi_o(z) \overset{\rm def}= \frac{\Psi(z)-\Psi(-z)}{2}.
\end{equation}

\paragraph{ }
Denote $\Phi(z)\colonequals F(z;\rho)$. The following identities are true
\begin{equation}\label{eq:85}
\begin{gathered}
\Phi_e'(z) = \frac 12 (\Phi'(z) - \Phi'(-z)) =\frac 12 (\Phi(\rho z) - \Phi(-\rho z))
= \Phi_o(\rho z),
\\
\Phi_o'(z) = \frac 12 (\Phi'(z) + \Phi'(-z)) =\frac 12 (\Phi(\rho z) + \Phi(-\rho z))
= \Phi_e(\rho z).
\end{gathered}
\end{equation}

Since $F(z;1)=e^{z}$, we obtain that $\Phi_e(z)=\cosh z$ and $\Phi_o(z)=\sinh z$ for $\rho=1$. In particular, $\Phi_e(iz)$ and $\Phi_o(iz)$ have simple real interlacing zeros in this case.

Further in this subsection we consider zeros of~$\Phi_e$ and~$\Phi_o$ for $0<\rho<1$. Denote
\begin{equation}\label{eq:chi}
\chi_e(z)\colonequals\Phi_e(\sqrt{z}) = \sum_{k=0}^\infty\frac{\rho^{n(2n-1)}z^n}{(2n)!},
\quad
\chi_o(z)\colonequals\dfrac{\Phi_o(\sqrt{z})}{\sqrt{z}} = \sum_{k=0}^\infty\frac{\rho^{n(2n+1)}z^n}{(2n+1)!}.
\end{equation}
\paragraph{ }
On one hand, $\Phi(z)$ is a stable real entire function of genus $0$. Therefore, $\Phi_e(z)$ and
$\Phi_o(z)$ have no common zeros. Indeed,
\begin{equation*}
\Phi_e(z_0)=\Phi_o(z_0)=0 \implies
\Phi_e(-z_0)=\Phi_o(-z_0)=0\implies
\Phi(z_0)=\Phi(-z_0)=0
\end{equation*}
is contradicting the fact that $\Phi(z)$ is stable (i.e. has no zeros in the closed complex
plane). Moreover, $\Phi_e(iz)$ and $\Phi_o(iz)$ have simple real interlacing zeros
(see~\cite{ChebMei}). That is, denoting the zeros of $\chi_o(z)$ and $\chi_e(z)$ as
$\{-\alpha_\nu\}_{\nu=1}^{\infty}$ and $\{-\beta_\nu\}_{\nu=1}^{\infty}$, respectively, we
obtain
\[
\chi_e(z) = \prod_{\nu=1}^\infty\left(1+\frac{z}{\beta_\nu}\right)
\quad\text{and}\quad
\chi_o(z) = \prod_{\nu=1}^\infty\left(1+\frac{z}{\alpha_\nu}\right),
\quad\text{where}\quad
0<\beta_1<\alpha_1<\beta_2<\alpha_2<\dots;
\]
here we took into account that $\Phi_o(0)=0$ to arrange the zeros.

\paragraph{ }
On the other hand, the equalities~\eqref{eq:85} give us
\[
\tfrac{d}{dz}(\Phi_e(iz)) = i\Phi_o(\rho iz),
\quad
\tfrac{d}{dz}(i\Phi_o(iz)) = -\Phi_e(\rho iz).
\]
Observe that $\Phi_e(iz)$ has the simple real zeros $\{\pm\sqrt{\beta_\nu}\}_{\nu=1}^\infty$ and
$\Phi_o(iz)$ has the simple real zeros $\{0\}\cup\{\pm\sqrt{\alpha_\nu}\}_{\nu=1}^\infty$.
According to Rolle's theorem ($\Phi_e(iz)$ and $i\Phi_o(iz)$ are real functions), we obtain that
\[
0<\beta_1<\rho^{-2}\alpha_1
<\beta_2<\rho^{-2}\alpha_2<\dots
\quad\text{and}\quad
0<\rho^{-2}\beta_1<\alpha_1
<\rho^{-2}\beta_2<\alpha_2<\dots
\]
which can be combined as follows
\begin{equation}\label{eq:10}
0<\beta_1<\rho^{-2}\beta_1
<\alpha_1<\rho^{-2}\alpha_1
<\beta_2<\rho^{-2}\beta_2
<\alpha_2<\rho^{-2}\alpha_2<\dots.
\end{equation}

\subsection{The conjecture for negative~$q$: the first approach}
\label{sec:conj-negat-q}

Here we assume $\rho>0$ and $q=-\rho$. To prove the conjecture the technique same as
in~\cite{Liu} can be applied, see Section~\ref{sec:conjecture-real-q} in the appendix. However,
since we already have obtained properties of $\Phi_e(z)$ and $\Phi_o(z)$ in the previous
subsection, this can be done much shorter. Note, that this case corresponds to the so-called
\emph{self-interlacing} polynomials of the paper\cite{Tyaglov}.

The formula~\eqref{eq:81a} gives (or, alternatively, this expression can be checked by the
direct substitution into~\eqref{eq:1})
\[F(z;-\rho)=\frac{1 - i}2 \Phi(i z) + \frac{1 + i}2 \Phi(-iz) =\Phi_e(i z) - i \Phi_o(iz).\]
The functions ${\Phi_e(i z)}$ and ${-i\Phi_o(iz)}$ are both real over the real line, and have
simple real interlacing zeros (see Subsection~\ref{sec:properties-even-odd})

Observe that the condition $F(z,-\rho)=0$ is equivalent to $\Phi_e(i z) = i \Phi_o(iz)$,
implying
\[R(z)\colonequals\frac{\Phi_e(i z)}{i\Phi_o(iz)}=1\]
The function $R$ has only real interlacing zeros and poles, so it maps the upper half of the
complex plane into itself or into the lower half of the complex plane. Therefore, if a complex
number~$z_0$ satisfies $R(z_0)=1\in\mathbb{R}$, then it must be real. Moreover, the derivative
$R'(z)$ is monotonic over the real line (see, \emph{e.g.}~\cite{ChebMei}). It is positive since
\[
\sign R'(z)
= \sign \big(i\Phi_o(i\rho 0)i\Phi_o(i0) + \Phi_e(i\rho 0) \Phi_e(i0)\big)
= 1.
\]
Thus, $R(z)=1$ once and only once in between consequent zeros and poles of $R(z)$ (each second
interval). Taking into account that $R(0)=\infty$ and $R(-z)=-R(z)$, we deduce that the first by
modulus solution to the equation~$R(z)=1$ (and hence to $F(z)=0$) is negative, then the next is positive,
then again negative, and so on.

\subsection{The case $q^2=-\rho^2$ for $\rho>0$}\label{sec:case-q2=-rho2}
\paragraph{ }
Take $\mu$ so that $\mu^2\rho=q$, in particular $\mu^4=-1$. In fact, this case is covered by
Corollary~\ref{cr:ca2}, however it is expository to apply Theorem~\ref{th:ca1} itself. The
formula~\eqref{eq:82} states that
\[
    F(z;\mu^2\rho) =\Phi_e(\mu z) + \mu\Phi_o(\overline\mu z).
\]
We can check this expression directly by substituting into the problem~\eqref{eq:1}: 
\begin{multline}\nonumber
    F'(z;q)
    = \mu\Phi_o(\rho\mu z) + \Phi_e(\rho\overline\mu z)
    = \Phi_e(-\rho\mu^3 z) + \mu\Phi_o(\rho\mu z)=\\
    \Phi_e(\pm i\theta\mu^3 \rho z) + \mu\Phi_o(\mu^{1-2}\theta\rho z))
    = \Phi_e(i\mu^3 (\theta\rho z)) + \mu\Phi_o(\overline\mu(\theta\rho z)) =F(qz;q)
\end{multline}
\[
\text{and}\quad
\Big[\Phi_e(\mu z) + \mu\Phi_o(\overline\mu z)\Big]_{z=0}=\Phi_e(0)=1.
\]

\paragraph{ }
The function
\[
\tilde{F}(z)\colonequals F(\overline{\mu} z,\mu^2\rho)=\Phi_e(z) + \mu\Phi_o(\mu^{-2}z) = \Phi_e(z) + \mu\Phi_o(\mu^{-2}z)
\]
has the form
\[
\tilde{F}(z)
= \Phi_e(z) + \mu^{-1}z\frac{\Phi_o(\mu^{-2}z)}{\mu^{-2}z}
= \Phi_e(z) + \overline{\mu} z\frac{\Phi_o(iz)}{iz} = \chi_e(z^2) + \overline{\mu} z\chi_o(-z^2)
\]
according to~\eqref{eq:chi}. Therefore, for $c=\overline{\mu}$ it satisfies one of the
cases~\ref{distr1} and \ref{distr2} from Theorem~\ref{th:ca1}. This implies that zeros of
$F(z)$ are located in such way, that Conjecture~\ref{conjecture-B} holds true.

You can find an illustration on the page~\pageref{fig:4}.

\newpage
\appendix
\section{Plots}\label{sec:graphics}

\begin{figure}[H]\label{fig:4}
    \includegraphics[scale=.9]{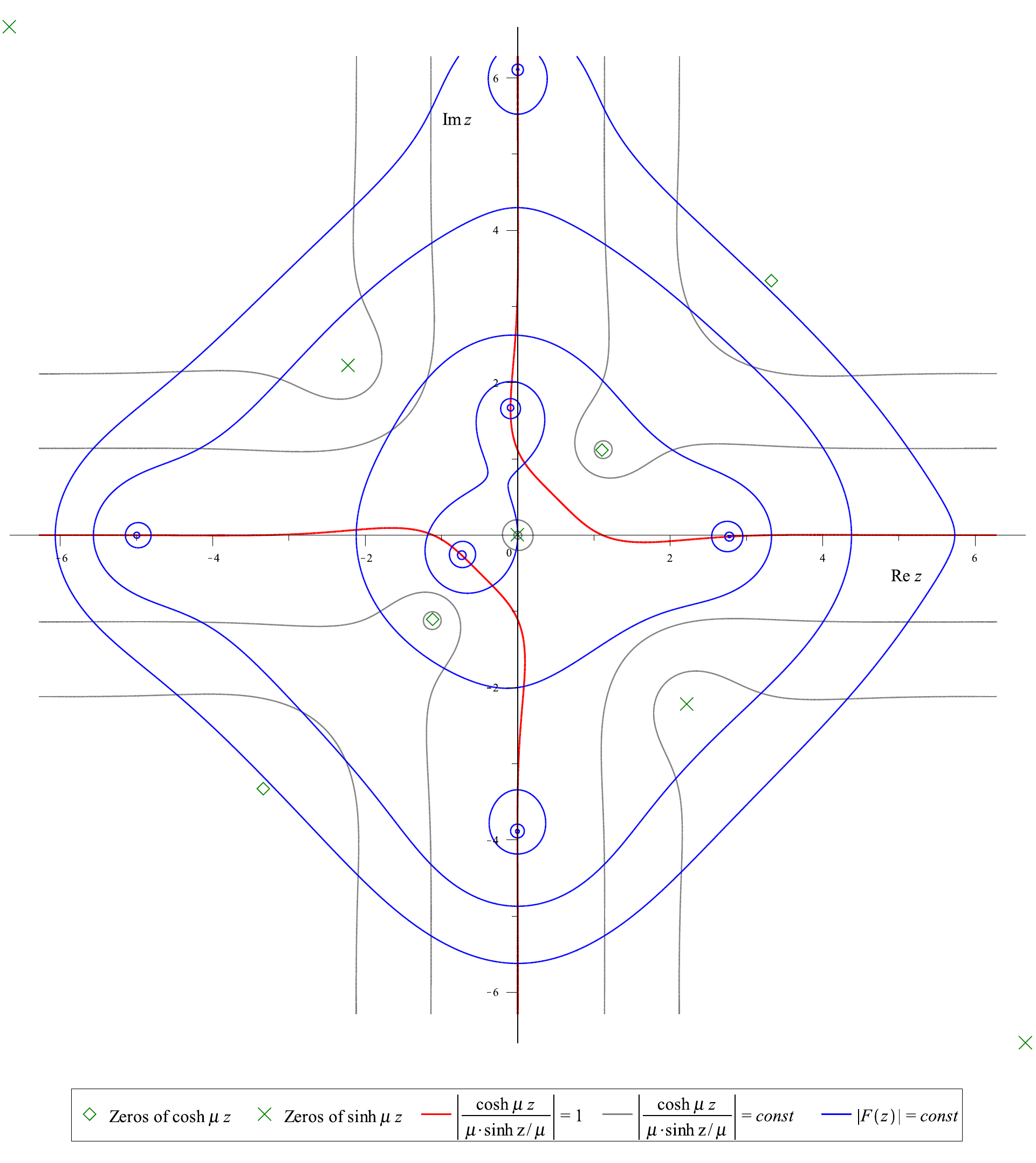}
    \caption{Illustration to the subsection~\ref{sec:case-q2=-rho2}.}
    \vskip .8em
    Localization of zeros of the function $F(z)=\cosh\mu z + \mu\cdot\sinh\dfrac{z}{\mu}$, where
    $\mu^2=i$.
    
    The zeros of $F(z)$ are inside of blue level lines; they all fall onto the red curve.
\end{figure}
\newpage

    \includegraphics[scale=.9]{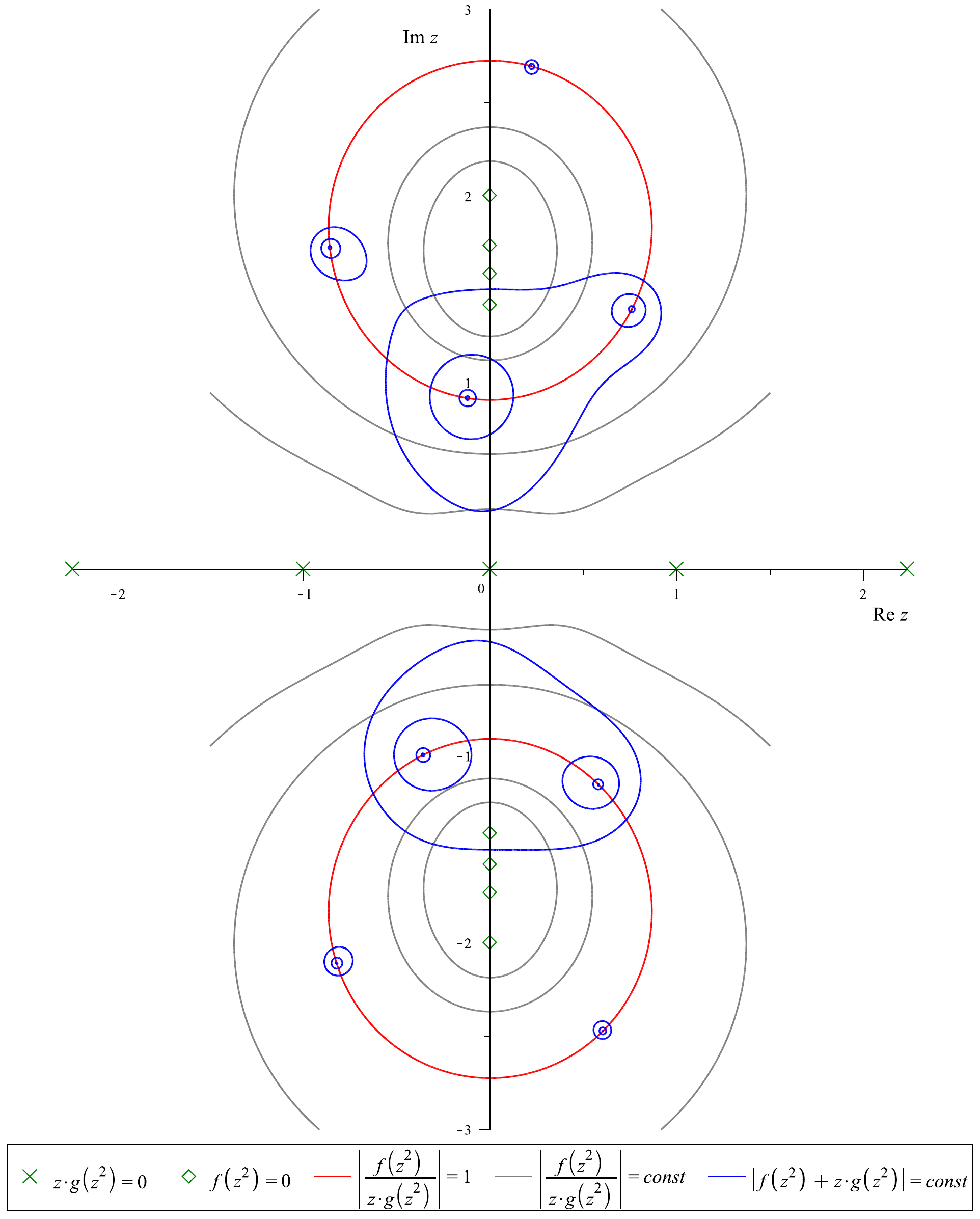}
\begin{figure}[H]\label{fig:3}
    \caption{Illustration to the proof of Theorem~\ref{th:ca1}.}
    \vskip .3em
    The case $f(z)=(z+2)\cdot(z+2.5)\cdot(z+3)\cdot(z+4)$ and $g(z)=(1+i)\cdot(z-1)\cdot(z-5)$,
    that is
    \[
    F(z)= f(z^2)+z\cdot g(z^2) = z^8 + 11.5z^6 + (1+i)z^5 +48.5 z^4 - 6(1+i)z^3 +
    89z^2 + 5(1+i) z + 60
    \]
\end{figure}

\section{Relation between different solutions to the problem~\eqref{eq:1}}
\subsection{The case~$q^n=\rho^n$ for any natural $n$}
Here in the case $q^n=\rho^n$, $n\in\mathbb{N}$, we obtain a relation between functions $F(z;q)$
and $F(z;\rho)$ of the complex variable~$z$ given by the equation~\eqref{eq:0}. The analogous
construction for the particular case of $|q|=1$ was implemented in~\cite{Alander,Valiron}.

\paragraph{ }
Assume that $q=\rho\mu^2$, where $\mu^{n}=-1$ and $|\rho|=|q|$. Note that whenever~$q^n=\rho^n$,
there exist one (if $n$ is odd) or two (if $n$ is even) corresponding values of~$\mu$. Moreover,
assume that $\mu^2$ is an $n$th \textbf{primitive} root of unity; in other words we choose a
minimal appropriate value of~$n$.

For the sake of brevity, denote $\Phi(z)\colonequals F(z;\rho)$; the function~$\Phi(z)$ is another
(``basic'') solution. We are going to express~$F(z;q)$ as follows
\begin{equation}\label{eq:6}
F(z;q)=\sum_{k=1}^{n}C_k\Phi(\mu^{2k+n-1}z).
\end{equation}
Since the Cauchy problem~\eqref{eq:1} is uniquely solvable, it is enough to find the
constants~$C_1, C_2, \dots, C_{n}$ such that~\eqref{eq:1} is satisfied by~\eqref{eq:6}.

\paragraph{ }
Substituting the representation~\eqref{eq:6} into the equation~\eqref{eq:1} gives
\begin{multline}\nonumber
    \sum_{k=1}^{n}C_k\Phi\big((\rho\mu^2)\mu^{2k+n-1}z\big)
    =\sum_{k=1}^{n}C_k\mu^{2k+n-1}\Phi'(\mu^{2k+n-1}z)\\
    =\sum_{k=1}^{n}C_k\mu^{2k+n-1}\Phi(\rho\mu^{2k+n-1}z)
    =\sum_{k=1}^{n}C_k\mu^{2k+n-1}\Phi\big((\rho\mu^2)\mu^{2(k-1)+n-1}z\big)\\
    =\sum_{k=1}^{n-1}C_{k+1}\mu^{2k+n+1}\Phi\big((\rho\mu^2)\mu^{2k+n-1}z\big) +
    C_{1}\mu^{n+1}\Phi\big((\rho\mu^2)\mu^{2n+n-1}z\big).
\end{multline}

\paragraph{ }
Observe from the last equality, that to satisfy~\eqref{eq:1} it is sufficient to pick up the
constants
\[
C_k  = C_{k+1}\mu^{2k+n+1}=-C_{k+1}\mu^{2k+1},\qquad k=1,2,\dots,n-1,\quad C_n=-C_{1}\mu,
\]
so that the initial condition $F(0;q)=1$ is true. So for $k=1,2,\dots,n$ we obtain
\[
C_k=-C_{k-1}\mu^{-(2k-1)}=\dots=(-1)^{k-1}C_1\mu^{-\sum_{j=2}^{k} (2j-1)}
=C_n(-1)^{k}\mu^{-\sum_{j=1}^{k} (2j-1)} =C_n(-1)^{k}\mu^{-k^2}.
\]

The last identity leaves the constant~$C_n$ unfixed, such that
\[
F(z;q)=C_n\sum_{k=1}^{n}(-1)^k\mu^{-k^2}\Phi(-\mu^{2k-1}z)
\]
Therefore, whenever
\[\sum_{k=1}^{n}(-1)^k\mu^{-k^2}\Phi(0) = \sum_{k=1}^{n}(-\mu)^{-k^2}\ne 0,\]
the initial condition can be satisfied by setting $C_n$ to be equal to
$\big(\sum_{k=1}^{n}(-\mu)^{-k^2}\big)^{-1}$. This leads us to the following representation of the
solution:
\begin{equation}
    \label{eq:77}
    \begin{gathered}
        F(z;q)=\nu_\mu\sum_{k=1}^{n}(-\mu)^{-k^2}\Phi(-\mu^{2k-1}z),\\
        \text{where}\quad \nu_\mu\colonequals\dfrac{1}{\sum_{k=1}^{n}(-\mu)^{-k^2}}.
    \end{gathered}
\end{equation}
The further steps \ref{step:n_odd} and \ref{step:n_even} are devoted to the calculation of the
coefficient~$\nu_\mu$. In particular, we show that
\[
\sum_{k=1}^{n}(-\mu)^{-k^2}\ne 0\quad\text{and}\quad|\nu_\mu|=\frac{1}{\sqrt{n}}.
\]

\paragraph{ }\label{step:n_odd}
Let $n=2m+1$ for some natural $m$. Observe that since $\mu^n=-1$, there exists an integer~$p$,
$-m\le p\le m$, such that
\[\mu=e^{(2p+1)\frac{\pi i}{2m+1}}.\]
Then
\[
\sum_{k=1}^{n}(-\mu)^{-k^2} = \sum_{k=1}^{n}\left(e^{(2p+1-n)\frac{\pi i}{n}}\right)^{-k^2}
    = \sum_{k=1}^{2m+1}e^{(m-p)\frac{2\pi i}{2m+1}k^2} = S(m-p,2m+1),
\]
where~$S(\makebox[.8em]{$\cdot$},\makebox[.8em]{$\cdot$})$ denotes the quadratic Gaussian sum
\[
S(m,n)\colonequals \sum_{k=0}^{n-1}e^{m\frac{2\pi i}{n}k^2} = \sum_{k=1}^{n}e^{m\frac{2\pi i}{n}k^2}.
\]

Since $\mu^2$ is an $n$-th primitive root of $1$, we know that $n-(2p+1) = 2(m-p)$ and $n=2m+1$
are coprime, which is equivalent to $\gcd(m-p,2m+1)=1$. Therefore, we obtain (see,
\emph{e.g.}~\cite[p.~26]{BRK_GS}):
\[
\sum_{k=1}^{n}(-\mu)^{-k^2}\in \left\{\sqrt{n},-\sqrt{n},i\sqrt{n},-i\sqrt{n}\right\},
\]
that is
\[
\nu_\mu^4 = n^{-2}\quad\text{for every}\quad n=1,3,5,\dots.
\]

\paragraph{ }\label{step:n_even}
Let $n$ be an even natural number. Then
\[
\sum_{k=n+1}^{2n}(-\mu)^{-k^2} = \sum_{k=1}^{n}(\mu^{n+1})^{-n^2-2kn}(-\mu)^{-k^2}
 = \sum_{k=1}^{n}1^{n+1}(-\mu)^{-k^2}
 = \sum_{k=1}^{n}(-\mu)^{-k^2}
\]
Using the representation
\[
\mu=e^{p\frac{\pi i}{n}},\quad\text{where}\quad 1\le p< n\quad\text{and}\quad\gcd(n-p,n)=1,
\]
we obtain
\[
\sum_{k=1}^{n}(-\mu)^{-k^2}=\frac 12 \sum_{k=1}^{2n}(-\mu)^{-k^2}
  =\frac 12 \sum_{k=1}^{2n}\left(e^{\frac{p-n}{n}\pi i}\right)^{-k^2}
  =\frac 12 \sum_{k=1}^{2n}e^{(n-p)\frac{2\pi i}{2n}k^2}
  =\frac 12 S(n-p,2n).
\]
Since $2n \equiv 0 \pmod 4$, the quadratic Gaussian sum $S(n-p,2n)$ has the following value (see,
\emph{e.g.}~\cite[p.~27]{BRK_GS}):
\[
S(n-p,2n)=\pm (1+i^{n-p})\sqrt{2n}.
\]

Note that since $n-p$ is odd, $i^{n-p}=\pm i$. Therefore, for even values of~$n$ we have
\[
\nu_\mu=1/\sum_{k=1}^{n}(-\mu)^{-k^2}\in
\left\{\frac{1+i}{\sqrt{2n}},\frac{-1-i}{\sqrt{2n}},
    \frac{1-i}{\sqrt{2n}},\frac{-1+i}{\sqrt{2n}}\right\},\quad \nu_\mu^4=-n^{-2}.
\]
\subsection{The case~$q^n=\rho^n$, where $n$ is even}
\paragraph{ }
For the even values of $n=2m$ we can express~\eqref{eq:77} in the equivalent form
\begin{multline}
    F(z;q)=\nu_\mu\sum_{k=1}^{m}\mu^{2km-k^2}\Phi(\mu^{2m+2k-1}z)+\nu_\mu\sum_{k=m+1}^{2m}\mu^{2km-k^2}\Phi(\mu^{2m+2k-1}z)\\
    =\nu_\mu\sum_{k=m+1}^{2m}\mu^{2(k-m)m-(k-m)^2}\Phi(\mu^{2m+2(k-m)-1}z)
    +\nu_\mu\sum_{k=1}^{m}\mu^{2(k+m)m-(k+m)^2}\Phi(\mu^{2m+2(k+m)-1}z)\\
    =\nu_\mu\sum_{k=m+1}^{2m}\mu^{-k^2+m^2+4m(k-m)}\Phi(\mu^{2k-1}z)+\nu_\mu\sum_{k=1}^{m}\mu^{-k^2+m^2}\Phi(\mu^{2k-1}z)\\
    =\nu_\mu\mu^{m^2}\sum_{k=1}^{2m}\mu^{-k^2}\Phi(\mu^{2k-1}z)
    =\tilde{\nu}_\mu\sum_{k=1}^{2m}\mu^{-k^2}\Phi(\mu^{2k-1}z),
\end{multline}
where $\tilde{\nu}_\mu\colonequals\mu^{m^2}\nu_\mu$.

\paragraph{ }
For $m=1$ we have $\mu=\pm i$,
\[
\tilde{\nu}_\mu^{-1}F(z;-\rho)
  =\mu^{-1}\Phi(\mu^{2-1}z)+\mu^{-4}\Phi(\mu^{4-1}z)
  =\overline\mu\Phi(\mu z)+\Phi(\overline\mu z)
  =\mp i \Phi(\pm i z) + \Phi(\mp i z)
\]
\begin{equation*}
\text{and}\quad \tilde{\nu}_2 = \tilde{\nu}_\mu=1/(1\mp i) = (1\pm i)/2 \qquad\Big[\tilde{\nu}_2 = (1+\mu)/2\Big];
\end{equation*}
so
\begin{equation*}
    F(z;-\rho) =\frac{(1\pm i)(\mp i)}2 \Phi(\pm i z) + \frac{1\pm i}2 \Phi(\mp i z)
    =\frac{1 - i}2 \Phi(i z) + \frac{1 + i}2 \Phi(-iz).
\end{equation*}
Applying of the notation~\eqref{eq:odd_even_def} gives
\begin{equation}\label{eq:81a}
    F(z;-\rho)
    =\frac{1}2 (\Phi(i z) + \Phi(-i z)) + \frac{i}2 (\Phi(iz) - \Phi(-iz))
    = \Phi_e(iz)+i\Phi_o(iz).
\end{equation}

\paragraph{ }
Analogously, for $m=2$ we have $\mu^4=-1$,
\begin{multline}\nonumber
    \tilde{\nu}_\mu^{-1}F(z;\mu^2\rho)
    =\mu^{-1}\Phi(\mu^{2-1}z)+\mu^{-4}\Phi(\mu^{4-1}z)+\mu^{-9}\Phi(\mu^{6-1}z)+\mu^{-16}\Phi(\mu^{8-1}z)\\
    =\overline\mu\Phi(\mu z)-\Phi(-\overline\mu z) +\overline\mu \Phi(- \mu z) + \Phi(\overline\mu z),
\end{multline}
\begin{equation*}
\tilde{\nu}_4=\tilde{\nu}_\mu=1/(\overline\mu + \overline\mu) = \mu/2,
\end{equation*}
so
\begin{equation}\label{eq:82}
    F(z;\mu^2\rho) =\tfrac 12 \big(\Phi(\mu z) + \Phi(- \mu z) + \mu\Phi(\overline\mu z) -
    \mu\Phi(-\overline\mu z) \big) =\Phi_e(\mu z) + \mu\Phi_o(\overline\mu z).
\end{equation}
where we applied the notation~\eqref{eq:odd_even_def}.


\paragraph{ }
Furthermore, if $m=4l$, then $\mu^{m^2}=i^{4l}=1$, so $\nu_\mu=\tilde{\nu}_\mu$ and
\[
F(-z;q)=\nu_\mu\sum_{k=1}^{8l}(-1)^k\mu^{-k^2}\Phi(\mu^{2k-1}z)
=\tilde{\nu}_\mu\sum_{k=1}^{8l}(-1)^k\mu^{-k^2}\Phi(\mu^{2k-1}z).
\]
Thus, we obtain $F(z;q)=F_e(z;q)+F_o(z;q)$, where
\begin{multline}\nonumber
    F_e(z;q) = \tilde{\nu}_\mu\sum_{k=1}^{4l}\mu^{-4k^2}\Phi(\mu^{4k-1}z) =
    \tilde{\nu}_\mu\sum_{k=1}^{2l}\left(\mu^{-4k^2}\Phi(\mu^{4k-1}z) +
        \mu^{-4(k+2l)^2}\Phi(-\mu^{4k-1}z)\right)\\
    = \tilde{\nu}_\mu\sum_{k=1}^{2l}\left(\mu^{-4k^2}\Phi(\mu^{4k-1}z) +
        \mu^{-4k^2-16l(k+l)}\Phi(-\mu^{4k-1}z)\right)
    = \tilde{\nu}_\mu\sum_{k=1}^{2l}\mu^{-4k^2}\Phi_e(\mu^{4k-1}z) ,
\end{multline}
\begin{multline}\nonumber
    F_o(z;q)= \tilde{\nu}_\mu\sum_{k=1}^{4l}\mu^{-(2k-1)^2}\Phi(\mu^{4k-3}z) =
    \tilde{\nu}_\mu\sum_{k=1}^{2l}\left(\mu^{-(2k-1)^2}\Phi(\mu^{4k-3}z) +
        \mu^{-(2k-1+4l)^2}\Phi(-\mu^{4k-3}z)\right)\\
    = \tilde{\nu}_\mu\sum_{k=1}^{2l}\left(\mu^{-(2k-1)^2}\Phi(\mu^{4k-3}z) +
        \mu^{-(2k-1)^2-16l(l+k)+8l}\Phi(-\mu^{4k-3}z)\right) \\
    = \tilde{\nu}_\mu\sum_{k=1}^{2l}\mu^{-(2k-1)^2}\left(\Phi(\mu^{4k-3}z) - \Phi(-\mu^{4k-3}z)\right) 
    = \tilde{\nu}_\mu\sum_{k=1}^{2l}\mu^{-(2k-1)^2}\Phi_o(\mu^{4k-3}z).
\end{multline}

\section{Conjecture~\ref{conjecture-B} for negative~$q$: the second approach}
\label{sec:conjecture-real-q}

\paragraph{ }
Since for $q=-1$ the series~\eqref{eq:0} has the form
\[
F(z;-1) = \sum_{k=0}^\infty \frac {(-1)^{\left\lfloor\frac{k}{2}\right\rfloor}}{k!}z^k = \cos z + \sin z,
\]
Conjecture~\ref{conjecture-B} holds true for~$q=-1$.

\paragraph{ }
Furthermore, let $-1<q<0$. We can apply the same approach as for $0<q<1$ in~\cite{Liu}. The
function~$F(z;-1)$ is of the Laguerre-P\'olya class since it has real only zeros and of genus
$1$. The function $G(z)\colonequals|q|^{\frac{1}{2}z(z-1)}=e^{\frac{1}{2}\ln|q|z(z-1)}$ is of the same
class since it has no zeros and $\frac{1}{2}\ln|q|<0$. Thus, by the P\'olya Theorem (see,
\emph{e.g.}~\cite[p.~43,~Satz~10.2]{Obreschkoff}) the function
\[
F(z)=F(z;q)
= \sum_{k=0}^\infty \frac {q^{\frac{k(k-1)}2}z^k}{k!}
= \sum_{k=0}^\infty \frac {(-1)^{\left\lfloor\frac{k}{2}\right\rfloor}G(k)}{k!}z^k,
\]
of genus $0$ is also of the Laguerre-P\'olya class, and hence can have real roots only. Let
$\{\lambda_k\}_{k=0}^\infty$,
\[ 0<|\lambda_0|\le|\lambda_1|\le|\lambda_2|\le\cdots, \]
be the zeros of $F(z)$.

 \paragraph{ }
Note, that $F(z)>0$ for $z\in(-|\lambda_0|,|\lambda_0|)$. According to \eqref{eq:1},
$F'(z)=F(qz)>0$ on the same interval, and hence $F'(\lambda_0)>0$, i.e. the zero is simple. We
have $\lambda_0<0$. Indeed, if $\lambda_0>0$ then
\[
F(\lambda_0)=F(0)+\int_0^{\lambda_0} F'(z)\,dz=1+\int_0^{\lambda_0} F(qz)\,dz>0
\text{ since }
qz\in(-|\lambda_0|,|\lambda_0|),
\]
which is impossible. In addition we obtain~$F(-\lambda_0)>0$, such that
$|\lambda_1|>|\lambda_0|$.

\paragraph{ }
Now, $\lambda_1>0$. Indeed, if $\lambda_1<0$ then $F(z)>0$ for $z\in(\lambda_0,|\lambda_1|)$.
Since $qz\in(0,|\lambda_1|)$ and $\lambda_1<\lambda_0<0$, that leads to the inequality
\[
F(\lambda_1)=\int_{\lambda_0}^{\lambda_1} F(qz)\,dz<0,
\]
contradicting to $F(\lambda_1)=0$. Furthermore, it implies~$F(-\lambda_1)<0$. For any
$z\in(0,\lambda_0/q]$,
\[
F(z)=F(0)+\int_0^z F(qt)\,dt > 0,
\]
therefore $\lambda_1>\lambda_0/q$. In particular, $F'(\lambda_1)=F(q\lambda_1)\ne 0$ provided by
$-\lambda_1<q\lambda_1<\lambda_0$. Since the zeros $\lambda_0$ and $\lambda_1$ are simple, the
function changes it sign to the opposite at those points.

\paragraph{ }
Denote $\lambda_{-1}=\lambda_{-2}=0$ and suppose we have proved that
\[
\lambda_{2j-2}<\lambda_{2j-4}<\dots<\lambda_{0}<0<\lambda_{1}<\lambda_{3}<\dots<\lambda_{2j-1},
\]
the zeros $\lambda_{0}$,\dots, $\lambda_{2j-1}$ are simple and $|\lambda_{k}|>|\lambda_{k-1}|/|q|$
for $k=1,\dots,2j-1$. Let us show that $\lambda_{2j}$ and $\lambda_{2j+1}$ satisfy these
relations.

\paragraph{ }\label{step:1}
Let $j=2m-1$ for $m=1,2,\dots$, then
\begin{equation}\label{eq:sign_of_F}
\begin{aligned}
    &F(z)<0&&\text{whenever}&&
    z\in(-|\lambda_{2j}|,\lambda_{2j-2})\cup(\lambda_{2j-1},|\lambda_{2j}|)&\text{and}\\
    &F(z)>0&&\text{whenever}&&
    z\in (\lambda_{2j-2},\lambda_{2j-4})\cup(\lambda_{2j-3},\lambda_{2j-1}).&
\end{aligned}
\end{equation}
If $\lambda_{2j}>0$, then we have the contradiction with~$F(\lambda_{2j})=0$:
\[
F(\lambda_{2j})=\int_{\lambda_{2j-1}}^{\lambda_{2j}} F(qz)\,dz<0, \text{ since }
qz\in(q\lambda_{2j},q\lambda_{2j-1})\subset(-|\lambda_{2j}|,\lambda_{2j-2}),
\]
hence $\lambda_{2j}<0$. Moreover, for any $z\in(\lambda_{2j-1}/q,\lambda_{2j-2})$ we have
$qz\in(\lambda_{2j-2}q,\lambda_{2j-1})$ and $F(qz)>0$. Consequently,
\[
F(z)=\int_{\lambda_{2j-2}}^z F(qt)\,dt > 0
\quad\text{whenever}\quad z\in[\lambda_{2j-1}/q,\lambda_{2j-2}).
\]
As a consequence, $|\lambda_{2j}|>|\lambda_{2j-1}|/|q|$. Since
$q\lambda_{2j}\in(\lambda_{2j-1},|\lambda_{2j}|)$ and $F'(\lambda_{2j})=F(q\lambda_{2j})<0$
according to~\eqref{eq:sign_of_F}, the zero $\lambda_{2j}$ is simple.

\paragraph{ }\label{step:2}
We have $F(z)<0$ for $z\in (\lambda_{2j-1},|\lambda_{2j+1}|)$ and $F(z)>0$ for
$z\in (-|\lambda_{2j+1}|,\lambda_{2j})$. If $\lambda_{2j+1}<0$, then we have the
contradiction
\[
F(\lambda_{2j+1})=\int_{\lambda_{2j}}^{\lambda_{2j+1}} F(qz)\,dz>0, \text{ since }
qz\in(q\lambda_{2j},q\lambda_{2j+1})\subset(\lambda_{2j-1},|\lambda_{2j+1}|),
\]
therefore $\lambda_{2j+1}>0$. Furthermore, for any
$z\in(\lambda_{2j-1},\lambda_{2j}/q)$ we have $qz\in(\lambda_{2j},\lambda_{2j-1}q)$, $F(qz)<0$
such that
\[
F(z)=\int_{\lambda_{2j-1}}^z F(qt)\,dt < 0
\quad\text{whenever}\quad z\in(\lambda_{2j-1},\lambda_{2j}/q].
\]
Therefore, the inequality~$|\lambda_{2j+1}|>|\lambda_{2j}|/|q|$ is valid. The zero
$\lambda_{2j+1}$ is simple provided by
\[
q\lambda_{2j+1}\in(-|\lambda_{2j+1}|,\lambda_{2j})
\overset{\eqref{eq:sign_of_F}\;}\implies
F'(\lambda_{2j+1})=F(q\lambda_{2j+1})>0.
\]

\paragraph{ }
If $j=2m$, we can repeat exactly the same steps \ref{step:1}, \ref{step:2} exchanging the sign
of $F(z)$, that is inverting the inequalities in~\eqref{eq:sign_of_F} as well as further where it
is required.

\end{document}